\documentclass[a4paper,12pt,reqno]{amsart}
\usepackage{amsmath}
\usepackage{amssymb}
\usepackage{amsthm}

\usepackage{amscd}
\usepackage{amsfonts}
\usepackage{setspace}
\usepackage{version}

\title{Almost sure large fluctuations of random multiplicative functions}
\author{Adam J Harper}
\address{Mathematics Institute, Zeeman Building, University of Warwick, Coventry CV4 7AL, England}
\email{A.Harper@warwick.ac.uk}
\date{31st December 2020}

\numberwithin{equation}{section}
\theoremstyle{plain}

\onehalfspace
\newcommand{\N}{\mathbb{N}}
\newcommand{\R}{\mathbb{R}}
\newcommand{\E}{\mathbb{E}}
\newcommand{\p}{\mathbb{P}}
\newcommand{\Z}{\mathbb{Z}}
\newcommand{\C}{\mathbb{C}}

\newtheorem{thmas1}{Theorem}
\newtheorem{thmas2}[thmas1]{Theorem}

\newtheorem{normapprox1}{Normal Approximation Result}
\newtheorem{normcomp1}{Normal Comparison Result}

\newtheorem{harman1}{Harmonic Analysis Result}
\newtheorem{numth1}{Number Theory Result}
\newtheorem{numth2}[numth1]{Number Theory Result}

\newtheorem{epres1}{Euler Product Result}
\newtheorem{epres2}[epres1]{Euler Product Result}

\newtheorem{chaos1}{Multiplicative Chaos Result}
\newtheorem{chaos2}[chaos1]{Multiplicative Chaos Result}
\newtheorem{chaos3}[chaos1]{Multiplicative Chaos Result}
\newtheorem{chaos4}[chaos1]{Multiplicative Chaos Result}

\newtheorem{prop1}{Proposition}
\newtheorem{prop2}[prop1]{Proposition}
\newtheorem{prop3}[prop1]{Proposition}

\begin{document}

\maketitle

\begin{abstract}
We prove that if $f(n)$ is a Steinhaus or Rademacher random multiplicative function, there almost surely exist arbitrarily large values of $x$ for which $|\sum_{n \leq x} f(n)| \geq \sqrt{x} (\log\log x)^{1/4+o(1)}$. This is the first such bound that grows faster than $\sqrt{x}$, answering a question of Hal\'asz and proving a conjecture of Erd\H{o}s. It is plausible that the exponent $1/4$ is sharp in this problem.

The proofs work by establishing a multivariate Gaussian approximation for the sums $\sum_{n \leq x} f(n)$ at a sequence of $x$, conditional on the behaviour of $f(p)$ for all except the largest primes $p$. The most difficult aspect is showing that the conditional covariances of the sums are usually small, so the corresponding Gaussians are usually roughly independent. These covariances are related to an Euler product (or multiplicative chaos) type integral twisted by additive characters, which we study using various tools including mean value estimates for Dirichlet polynomials, high mixed moment estimates for random Euler products, and barrier arguments with random walks.
\end{abstract}

% INTRODUCTION %%%%%%%%%%%%%%%%%%%%%%%%%%%%%
\section{Introduction}
There are at least two different models of random multiplicative function $f(n)$ that have attracted quite a lot of attention in number theory, analysis and probability. A {\em Steinhaus random multiplicative function} is obtained by letting $(f(p))_{p \; \text{prime}}$ be a sequence of independent Steinhaus random variables (i.e. distributed uniformly on the unit circle $\{|z|=1\}$), and then setting $f(n) := \prod_{p^{a} || n} f(p)^{a}$ for all $n \in \N$, where $p^a || n$ means that $p^a$ is the highest power of the prime $p$ that divides $n$. We obtain a {\em Rademacher random multiplicative function} by letting $(f(p))_{p \; \text{prime}}$ be independent Rademacher random variables (i.e. taking values $\pm 1$ with probability $1/2$ each), and setting $f(n) := \prod_{p |n} f(p)$ for all squarefree $n$, and $f(n) = 0$ when $n$ is not squarefree. Rademacher random multiplicative functions were introduced by Wintner~\cite{wintner} as a model for the M\"{o}bius function $\mu(n)$. Steinhaus random multiplicative functions arise as models for randomly chosen Dirichlet characters $\chi(n)$ or ``continuous characters'' $n^{-it}$, see e.g. section 2 of Granville and Soundararajan~\cite{gransoundlcs}. The introduction to~\cite{harperrmflow} contains a more extensive survey of the literature in this area.

\vspace{12pt}
A fundamental type of result in the classical study of sums of independent random variables is the Law of the Iterated Logarithm, which describes the almost sure size of the largest fluctuations of those sums. For example, if $(\epsilon_n)_{n=1}^{\infty}$ is a sequence of independent random variables taking values $\pm 1$ with probability $1/2$ each, then Khintchine's law of the iterated logarithm asserts that almost surely, we have
$$ \limsup_{x \rightarrow \infty} \frac{\sum_{n \leq x} \epsilon_n}{\sqrt{2x\log\log x}} = 1 \;\;\; \text{and} \;\;\; \liminf_{x \rightarrow \infty} \frac{\sum_{n \leq x} \epsilon_n}{\sqrt{2x\log\log x}} = -1 . $$
See e.g. Chapter 8 of Gut~\cite{gut} for much more discussion of results of this flavour. Note that the largest fluctuations as $x$ varies are significantly larger than the random fluctuations one expects at any fixed point $x$. Indeed, for large $x$ the sum $\sum_{n \leq x} \epsilon_n$ will have a roughly Gaussian distribution with mean zero and variance $x$, so will typically have size like a constant times $\sqrt{x}$ (rather than a constant times $\sqrt{x\log\log x}$). Thus the Law of the Iterated Logarithm encodes information, not only about the distribution at each fixed point, but about the interaction between the sums as $x$ varies.

Random multiplicative functions cannot be studied using classical theory for independent random variables, because the values $f(n)$ are not independent (e.g. $f(6) = f(2)f(3)$ is entirely determined by $f(2)$ and $f(3)$). However, much work has been done developing alternative methods to study $\sum_{n \leq x} f(n)$. Wintner~\cite{wintner} analysed the convergence of the random Dirichlet series and Euler product corresponding to $f(n)$, and deduced (in the Rademacher case) that for any fixed $\epsilon > 0$, one almost surely has $\sum_{n \leq x} f(n) = O(x^{1/2+\epsilon})$ and one almost surely does {\em not} have $\sum_{n \leq x} f(n) = O(x^{1/2-\epsilon})$. In unpublished work reported in his problems paper~\cite{erdos}, Erd\H{o}s proved that the sum is almost surely $O(\sqrt{x} \log^{A}x)$ and almost surely not $O(\sqrt{x} \log^{-B}x)$. In the 1980s, Hal\'asz~\cite{halasz} introduced several beautiful new ideas to make further progress. Thus he made systematic use of conditioning, coupled with hypercontractive inequalities (i.e. high moment bounds) applied to the randomness coming from the large primes, to prove that almost surely (in the Rademacher case) one has a bound $O(\sqrt{x} e^{A\sqrt{\log\log x \log\log\log x}})$. Using the classical connection between integral averages of $\sum_{n \leq x} f(n)$ and of the random Euler product corresponding to $f(n)$, coupled with an ingenious special argument for studying suprema of such products, he also proved that one almost surely does {\em not} have a bound $O(\sqrt{x} e^{-B\sqrt{\log\log x \log\log\log x}})$. (It is worth noting that the similar shape of Hal\'asz's upper and lower bounds is basically coincidental, as the methods involved are quite distinct.)

Most recently, Lau, Tenenbaum and Wu~\cite{tenenbaum} (see also Basquin's paper~\cite{basquin}) introduced a splitting device into Hal\'asz's upper bound argument. This significantly increased the efficiency of the hypercontractive inequalities, yielding an almost sure bound $O(\sqrt{x} (\log\log x)^{2+\epsilon})$. For the lower bound, Harper~\cite{harpergp} introduced Gaussian process tools to sharpen Hal\'asz's analysis of suprema of random Euler products, establishing (again in the Rademacher case) that almost surely $\sum_{n \leq x} f(n) \neq O(\sqrt{x}/(\log\log x)^{5/2 + \epsilon})$.

These results begin to resemble the Law of the Iterated Logarithm, but are still rather far from that in terms of precision. In particular, the lower bounds do not exclude the striking (and at first sight unlikely looking) possibility that one might almost surely have $\sum_{n \leq x} f(n) = O(\sqrt{x})$, or even $\sum_{n \leq x} f(n) = o(\sqrt{x})$. Hal\'asz posed this as a question in the Steinhaus case (see the Unsolved Problems in Montgomery's book~\cite{mont}), and it is given further credence by the recently proved fact that one does see better than squareroot cancellation in low moments of random multiplicative functions (e.g. $\E|\sum_{n \leq x} f(n)| \asymp \frac{\sqrt{x}}{(\log\log x)^{1/4}}$, see Harper's paper~\cite{harperrmflow}). On the other hand, Erd\H{o}s~\cite{erdos} conjectured (in the Rademacher case) that almost surely $\sum_{n \leq x} f(n) \neq O(\sqrt{x})$, whilst writing that ``As far as I know nobody has a plausible guess for the true order of magnitude of $|\sum_{n \leq x} f(n)|$''. In the Rademacher case, where $f(n)$ was originally intended as a model for the M\"{o}bius function $\mu(n)$, one might also have in mind Gonek's conjecture (see Ng's paper~\cite{ng}) that the largest fluctuations of $\sum_{n \leq x} \mu(n)$ should have size $\sqrt{x} (\log\log\log x)^{5/4}$. Existing results on $\sum_{n \leq x} f(n)$ are consistent with such behaviour.

\vspace{12pt}
In this paper we shall prove the following result, answering Hal\'asz's question in the negative and proving the conjecture of Erd\H{o}s.

\begin{thmas1}\label{thmglobalas}
Let $f(n)$ be a Steinhaus or Rademacher random multiplicative function. Then for any function $V(x)$ tending to infinity with $x$, there almost surely exist arbitrarily large values of $x$ for which
$$ \Bigl| \sum_{n \leq x} f(n) \Bigr| \geq \frac{\sqrt{x} (\log\log x)^{1/4}}{V(x)} . $$
\end{thmas1}

As we will explain, it now seems plausible that a bound of the form $\sqrt{x} (\log\log x)^{1/4 + o(1)}$ is sharp in this problem, and if so Theorem \ref{thmglobalas} is quite close to best possible.

Classically, to prove the ``lower bound'' part of a Law of the Iterated Logarithm (i.e. that large fluctuations almost surely occur) one would construct a sequence of {\em independent} events that force the sum under consideration to be large at different points $x$, such that the sum of their probabilities diverges, and apply the second Borel--Cantelli lemma to conclude that almost surely infinitely many of the events will occur. For example, one might consider events that $\sum_{x_{j-1} < n \leq x_j} \epsilon_n$ is large, for a suitable rapidly increasing sequence $x_j$. It seems difficult to apply this strategy to random multiplicative functions, because there is substantial dependence between the sums $\sum_{n \leq x} f(n)$ even for very widely separated $x$ (although it is possible to understand the nature of this dependence and remove some of it via conditioning, as we shall see). Instead we shall organise our arguments differently, with most of our work going into proving the following localised statement.

\begin{thmas2}\label{thmlocalas}
Let $f(n)$ be a Steinhaus or Rademacher random multiplicative function. Then uniformly for all large $X$ and all $1 \leq W \leq (1/30)\log\log\log X$, say, we have
$$ \max_{X^{8/7} \leq x \leq X^{4/3}} \frac{1}{\sqrt{x}} \Bigl| \sum_{n \leq x} f(n) \Bigr| \geq \frac{(\log\log X)^{1/4}}{e^{1.2W}} $$
with probability $\geq 1 - O(e^{-0.1W})$.
\end{thmas2}

\begin{proof}[Proof of Theorem \ref{thmglobalas}, assuming Theorem \ref{thmlocalas}]
We may assume that $V(x) \leq (\log\log x)^{1/25}$ for all large $x$, since the statement of Theorem \ref{thmglobalas} is stronger when $V(x)$ grows more slowly. Let us define $W(X) := \min\{\frac{\log V(x)}{1.2} - 1 : X^{8/7} \leq x \leq X^{4/3}\}$, noting that $W(X) \leq (1/30)\log\log\log X$ for all large $X$, and that since $V(x) \rightarrow \infty$ as $x \rightarrow \infty$ we also have $W(X) \rightarrow \infty$ as $X \rightarrow \infty$. Theorem \ref{thmlocalas} then implies that the probability of {\em not} finding a value $X^{8/7} \leq x \leq X^{4/3}$ for which $|\sum_{n \leq x} f(n)| \geq \sqrt{x} \frac{(\log\log X)^{1/4}}{e^{1.2W(X)}} \geq \frac{\sqrt{x} (\log\log x)^{1/4}}{V(x)}$ is $O(e^{-0.1W(X)})$. Summing this up along a suitably sparse sequence of $X$-values (depending on the size of $W(X)$), the sum will converge and so the first Borel--Cantelli lemma implies that almost surely, only finitely many of these ``failure'' events will occur.
\end{proof}

Note that the deduction of Theorem \ref{thmglobalas} here only required the {\em first} Borel--Cantelli lemma, which needs no independence assumptions.

\vspace{12pt}
To explain how Theorem \ref{thmlocalas} is proved, first recall how $\sum_{n \leq x} f(n)$ behaves for fixed $x$. As explored in the author's paper~\cite{harperrmflow}, the size of $\frac{1}{\sqrt{x}} \sum_{n \leq x} f(n)$ for $x$ around $X$ is closely related to the size of integrals like $\int_{-1/2}^{1/2} |F(1/2+it)|^2 dt$, where $F(s)$ is the random Euler product corresponding to $f(n)$ on numbers with all their prime factors at most $X$ (so $F(s) := \prod_{p \leq X} (1 - \frac{f(p)}{p^s})^{-1}$ in the Steinhaus case and $F(s) := \prod_{p \leq X} (1 + \frac{f(p)}{p^s})$ in the Rademacher case). One can investigate the distribution of the random integral $\int_{-1/2}^{1/2} |F(1/2+it)|^2 dt$ using ideas from the theory of {\em critical multiplicative chaos}, as we indicate later (and see section \ref{subsecmultchaos} below). In particular, one can compute moments of such integrals, and deduce (see Corollary 2 of Harper~\cite{harperrmflow}) that for all $\lambda \geq 2$, we have
$$ \p\Biggl(|\sum_{n \leq x} f(n)| \geq \lambda \frac{\sqrt{x}}{(\log\log x)^{1/4}} \Biggr) \ll \frac{\min\{\log\lambda, \sqrt{\log\log x}\}}{\lambda^2} . $$
Thus the sum will typically have size $\ll \frac{\sqrt{x}}{(\log\log x)^{1/4}}$, and to prove Theorem \ref{thmlocalas} we must use the maximum over $x$ to produce significantly larger values.

To achieve this, the first key idea is to exploit the randomness coming from $f(p)$ on ``large'' primes $p$ to a much greater extent than previously. More precisely, for all $X^{8/7} \leq x \leq X^{4/3}$ (say) we can use multiplicativity of $f(n)$ to write
$$ \frac{1}{\sqrt{x}} \sum_{n \leq x} f(n) = \frac{1}{\sqrt{x}} \sum_{X < p \leq x} f(p) \sum_{m \leq x/p} f(m) + \frac{1}{\sqrt{x}} \sum_{\substack{n \leq x, \\ n \; \text{is} \; X \; \text{smooth}}} f(n) . $$
Recall that a number is said to be $X$-smooth if all of its prime factors are $\leq X$. In previous work (e.g. in Hal\'asz's upper bound~\cite{halasz}, and the moment results of Harper~\cite{harperrmflow, harperrmfhigh}), decompositions like this were used along with conditioning on the behaviour of $f(p)$ for ``small'' $p$ (so the inner sums $\sum_{m \leq x/p} f(m)$ could be treated as fixed), and moment bounds for weighted sums like $\sum_{X < p \leq x} f(p) a_p$, to ultimately deduce moment or large deviation bounds for $\sum_{n \leq x} f(n)$. Here we instead prove a {\em multivariate Gaussian approximation} for (the real parts of) the sums $\frac{1}{\sqrt{x}} \sum_{X < p \leq x} f(p) \sum_{m \leq x/p} f(m)$, conditional on the values $(f(p))_{p \leq X}$, at a well spaced sequence of points $X^{8/7} \leq x \leq X^{4/3}$. In fact it isn't too difficult to obtain this much stronger information about the conditional distribution, as after conditioning the inner sums $\frac{1}{\sqrt{x}} \sum_{m \leq x/p} f(m)$ are fixed quantities with typical size $\approx \frac{1}{\sqrt{p}}$, so we just need a multivariate Gaussian approximation for the sum of many independent random variables $f(p)$ multiplied by various small coefficients. See Normal Approximation Result 1 in the Tools section (section \ref{sectools} below) for the Gaussian approximation theorem we use, which is applied in section \ref{subsecmostofproof}.

The key advantage of passing to multivariate Gaussian distributions is that to understand their behaviour, including existence of large values, one only needs (in principle) to understand their means, variances, and covariances. Here the means are always zero, irrespective of the values of $(f(p))_{p \leq X}$ on which we condition. The conditional variance at each $X^{8/7} \leq x \leq X^{4/3}$ may be shown, after a little calculation, to be $\gg \frac{1}{\log X} \int_{1}^{X^{1/7}} |\sum_{m \leq z} f(m)|^2 \frac{dz}{z^2}$. Note that one has {\em the same integral} for all such $x$, indicating some of the significant dependence between the corresponding sums extracted by conditioning on $(f(p))_{p \leq X}$. After applying a form of Parseval's identity, and some further calculation, this scaled integral may be lower bounded by something like $\frac{1}{\log X} \int_{1/3}^{1/2} |F(1/2+it)|^2 dt$, where $F(s)$ is the random Euler product. The size of this integral certainly does depend on the $(f(p))_{p \leq X}$, but one can show that with probability $1 - O(e^{-0.1W})$ those $f(p)$ will be such that $\frac{1}{\log X} \int_{1/3}^{1/2} |F(1/2+it)|^2 dt \geq \frac{1}{e^{2.2W} \sqrt{\log\log X}}$. These calculations are similar to some performed by the author~\cite{harperrmflow} when investigating moments of $\sum_{n \leq x} f(n)$, but in the details a fair amount of new work is required to obtain a lower bound that holds with probability close to 1 as opposed to positive probability. In particular, one needs a different treatment of the contribution to $F(s)$ from small primes, and one needs to work with the {\em variance of the integral} (with suitable barrier conditions imposed) as opposed to its cruder mean square. See Multiplicative Chaos Results 3 and 4 in section \ref{sectools}, which are applied in section \ref{subsecmostofproof}.

Now if the sums $\Re \frac{1}{\sqrt{x}} \sum_{X < p \leq x} f(p) \sum_{m \leq x/p} f(m)$, conditional on the values $(f(p))_{p \leq X}$, usually behaved roughly independently at the $\asymp \log X$ different sample points $x = X^{8/7} e^{2\pi r}$ (say) in the interval $[X^{8/7},X^{4/3}]$, then one would expect the maximum of those sums to behave like the maximum of $\asymp \log X$ independent Gaussians, each with mean zero and variance $\gg \frac{1}{e^{2.2W} \sqrt{\log\log X}}$. A standard calculation shows that the maximum of $n$ independent, mean zero, variance one Gaussians is $\approx \sqrt{2\log n}$ with high probability, so the maximum of our sums would be $\gg \sqrt{\log\log X} \cdot \sqrt{\frac{1}{e^{2.2W} \sqrt{\log\log X}}} \gg \frac{(\log\log X)^{1/4}}{e^{1.1W}}$ with high probability, as desired for Theorem \ref{thmlocalas}. The actual bound there is slightly weaker because one needs to account for various constants and other minor factors. And crucially, to assess whether multivariate Gaussians behave roughly independently one only needs to know that their covariances are sufficiently small compared with their variances. See Normal Comparison Result 1 in section \ref{sectools}, which we apply in section \ref{subsecmostofproof}. Thus the main technical work in proving Theorem \ref{thmlocalas} lies in bounding those covariances.

Similarly as above, one can show that the conditional covariance at two points $X^{8/7} \leq x, y \leq X^{4/3}$ is $\approx \frac{1}{\log X} \int_{-1/2}^{1/2} |F(1/2+it)|^2 (x/y)^{it} dt$. In fact it requires quite a lot of analytic work, using Perron's formula and Dirichlet polynomial estimates, to bring the conditional covariance into roughly this form. See section \ref{subseccovarmanip}. In particular, if one looks at points $x = X^{8/7} e^{2\pi r}$ and $y = X^{8/7} e^{2\pi u}$ then the covariance has the rough shape $\frac{1}{\log X} \int_{-1/2}^{1/2} |F(1/2+it)|^2 e^{2\pi i(r-u)t} dt$, and one wants to show that for most pairs $r \neq u$ this will be smaller than our lower bound $\frac{1}{e^{2.2W} \sqrt{\log\log X}}$ for the variance. To achieve this, an initial idea is to seek strong upper bounds (holding with high probability over $(f(p))_{p \leq X}$) for a quantity like $\sum_{u} \left| \frac{1}{\log X} \int_{-1/2}^{1/2} |F(1/2+it)|^2 e^{2\pi iut} dt  \right|^{2k}$. If one expands the $2k$-th power and sums over the $\asymp \log X$ values of $u$, one gets a big saving for all $t_1, ..., t_{2k}$ except where $t_1 + ... + t_k - t_{k+1} - ... - t_{2k}$ is very close (within distance $1/\log^{2/3}X$, say) to an integer. These bad tuples form a small subset of $[-1/2,1/2]^{2k}$, and to conclude one wants to know that the integral of the Euler products $\frac{1}{\log X} |F(1/2+it_j)|^2$ over this bad set will itself be small (i.e. much smaller than $(\frac{1}{e^{2.2W} \sqrt{\log\log X}})^{2k}$) with high probability.

But matters now become very subtle, because typically one finds significant contributions to $\frac{1}{\log X} \int_{-1/2}^{1/2} |F(1/2+it)|^2 dt$ from large values of $F(1/2+it)$ on small intervals of $t$ (much smaller than $1/\log^{0.99}X$, say), and a careful analysis of the structure of those contributions is required. To succeed we first show that with probability close to 1, we may discard the contribution from points $t$ where $|F(1/2+it)|$ {\em or any of its sufficiently long partial Euler products} are within a factor $\frac{1}{(\log\log X)^{1000}}$ of their mean square. Simultaneously considering partial products is heavily motivated by multiplicative chaos ideas, and the proof that we may restrict their sizes ultimately relies on two random walk calculations, one involving the maximum of the random walk and one connected with the so-called Ballot Theorem. See Multiplicative Chaos Results 1 and 2 in the Tools section (section \ref{sectools}), which are deployed for our main proof in section \ref{subsecstrongbarrier}.

With this modification, after expanding out $\sum_{u} \left| \frac{1}{\log X} \int_{-1/2}^{1/2} |F(1/2+it)|^2 e^{2\pi iut} dt  \right|^{2k}$ one can distinguish two kinds of contribution: from tuples $t_1, ..., t_{2k}$ where many consecutive $t_j$ are close together, and from tuples where many of the gaps between consecutive $t_j$ are fairly large. In both cases, to bound the contribution one first pulls out partial Euler products over suitable ranges of primes, and uses our restricted upper bound for those. (One pulls out more primes for points $t_j$ that are more clustered together, so that the products remaining are sufficiently decorrelated that taking expectations later will yield reasonable estimates. This again is a typical multiplicative chaos style of argument.) If many $t_j$ are close together, this process extracts many long partial products and so many saving factors $\frac{1}{(\log\log X)^{1000}}$, which suffice for the final bound. In contrast, if there are many large gaps then the condition that $t_1 + ... + t_k - t_{k+1} - ... - t_{2k}$ is within distance $1/\log^{2/3}X$ of an integer severely shrinks the corresponding set of $t_j$ (which it may {\em not} in the case where many $t_j$ are already close together), and this saving suffices for the final bound after taking expectations and applying Markov's inequality. To implement all of this we need, in particular, uniform bounds for very high moments of shifted Euler products (as explained in section \ref{subsecfinalexp}, one needs to take $k \asymp (\log\log X)/(\log\log\log X)$ in the above setup to deduce that the conditional covariance is small for sufficiently many pairs of $r,u$). See Euler Product Results 1 and 2 in section \ref{sectools}, which are deployed for our main proof in section \ref{subsecfinalexp}.

\vspace{12pt}
We finish with two remarks. Firstly, note that in the classical Law of the Iterated Logarithm one obtains fluctuations of size $\asymp \sqrt{x\log\log x}$ because the typical size of $\sum_{n \leq x} \epsilon_n$ is $\asymp \sqrt{x}$, and taking a maximum over well spaced points $x_j$ (one must at least have $x_{j+1}/x_{j} \geq \lambda > 1$ for the sums up to $x_j$ to be sufficiently independent of one another) multiplies this by a factor $\asymp \sqrt{\log\log x}$. In the random multiplicative case the typical size of $\sum_{n \leq x} f(n)$ is $\asymp \frac{\sqrt{x}}{(\log\log x)^{1/4}}$, and taking a maximum over similarly well spaced points again multiplies this by a factor $\approx \sqrt{\log\log x}$, as in Theorem \ref{thmlocalas}. This explains our earlier remark that a bound $\sqrt{x} (\log\log x)^{1/4 + o(1)}$ is plausibly sharp in Theorem \ref{thmglobalas}. However, if that is true it remains a formidable task to prove it. As will be seen in section \ref{subsecmostofproof}, for our lower bound we can discard the sum $\sum_{\substack{n \leq x, \\ n \; \text{is} \; X \; \text{smooth}}} f(n)$ fairly trivially (using the fact that there usually won't be many $x$ where it is remarkably large) and rely on the numbers with a prime factor $> X$ to produce large fluctuations, but for an upper bound one needs a good understanding of {\em both} sums for {\em all} $x$.

Secondly, one might ask whether it is possible to formulate, and hopefully prove, deterministic number theoretic statements corresponding to our Theorems. For Theorem \ref{thmglobalas} this isn't so clear, since the notion of something happening almost surely for a class of infinite sequences doesn't translate very naturally to finitary number theoretic objects. But for Theorem \ref{thmlocalas} the analogy is easier to draw: one should seek lower bounds for $\max_{X^{8/7} \leq x \leq X^{4/3}} \frac{|\sum_{n \leq x} \chi(n)|}{\sqrt{x}}$ that hold for ``most'' (e.g. for $99\%$ of) Dirichlet characters $\chi$ modulo a large prime $r$, say. Although the distribution of character sums $\sum_{n \leq x} \chi(n)$ has been extensively investigated (see~\cite{bobgoldgrankou, gransoundlcs}), this problem of producing large values by varying the endpoint in an interval (as opposed to varying the character $\chi$) seems relatively unexplored. If $X$ is small enough compared with $r$, one expects things to be well modelled by the Steinhaus case of Theorem \ref{thmlocalas}, and to have similar behaviour as there. On the other hand, we know that for most $\chi$ mod $r$ one has $\max_{x \in \N} |\sum_{n \leq x} \chi(n)| \ll \sqrt{r}$ (the proportion of $\chi$ for which this holds can be made arbitrarily close to 1 by increasing the implicit constant, see Bober, Goldmakher, Granville and Koukoulopoulos~\cite{bobgoldgrankou} for the state of the art on this problem). Thus the direct analogue of Theorem \ref{thmlocalas} cannot hold once $X^{8/7} \gg r/\sqrt{\log\log r}$, although one wouldn't expect it to since for such $X$, most of the range $X^{8/7} \leq x \leq X^{4/3}$ consists of values $x \geq 2r$ that contribute negligibly to the maximum due to periodicity of $\chi$ mod $r$. It would be interesting to determine whether a version of Theorem \ref{thmlocalas} does hold for character sums when $X \leq r^{3/4}$, say (so that $X^{4/3} \leq r$), and to investigate the transition in behaviour as $X$ becomes even larger.

\subsection{Notation and references}
We write $f(x) = O(g(x))$ and $f(x) \ll g(x)$, both of which mean that there exists $C$ such that $|f(x)| \leq Cg(x)$, for all $x$. We write $f(x) \asymp g(x)$ to mean that $g(x) \ll f(x) \ll g(x)$, in other words that $c g(x) \leq |f(x)| \leq C g(x)$ for some $c,C$, for all $x$.

The books of Gut~\cite{gut} and of Montgomery and Vaughan~\cite{mv} are excellent general references for probabilistic and number theoretic background for this paper.

% SECTION 2 %%%%%%%%%%%%%%%%%%%%%%%%%%%%%
\section{Tools}\label{sectools}

\subsection{Gaussian random variables}
As discussed in the Introduction, and as in many probabilistic problems, the basic structure of our proof will be to approximate our random variables of interest fairly precisely by Gaussian random variables (after performing a large amount of conditioning), at which point we can access some of the many special tools available for studying Gaussians. We present suitable results here, for use later.

For the approximation step, we need a multivariate central limit theorem with quantitative error terms. By now this is fairly standard in the setting we require, with the following result ultimately flowing from Stein's method of exchangeable pairs.

\begin{normapprox1}
Suppose that $m \geq 1$, and that $\mathcal{H}$ is a finite non-empty set. Suppose that for each $1 \leq i \leq m$ and $h \in \mathcal{H}$ we are given a deterministic coefficient $c(i,h) \in \C$. Finally, suppose that $(V_i)_{1 \leq i \leq m}$ is a sequence of independent, mean zero, complex valued random variables, and let $Y = (Y_h)_{h \in \mathcal{H}}$ be the $\#\mathcal{H}$-dimensional random vector with components $Y_h := \Re\left( \sum_{i=1}^{m} c(i,h) V_i \right)$.

If $Z = (Z_h)_{h \in \mathcal{H}}$ is a multivariate normal random vector with the same mean vector and covariance matrix as $Y$, then for any $u \in \R$ and any small $\eta > 0$ we have
\begin{eqnarray}
\p(\max_{h \in \mathcal{H}} Y_h \leq u ) & \leq & \p(\max_{h \in \mathcal{H}} Z_h \leq u + \eta ) + \nonumber \\
&& + O\left(\frac{1}{\eta^{2}} \sum_{g,h \in \mathcal{H}} \sqrt{\sum_{i=1}^{m} |c(i,g)|^2 |c(i,h)|^2 \E|V_i|^4 } + \frac{1}{\eta^3} \sum_{i=1}^{m} \E|V_i|^3 (\sum_{h \in \mathcal{H}} |c(i,h)|)^3 \right) . \nonumber
\end{eqnarray}
\end{normapprox1}

\begin{proof}[Proof of Normal Approximation Result 1]
This is a special case of Central Limit Theorem 1 from Appendix B of Harper~\cite{harperlcz}, which is itself an application of Theorem 2.1 of Reinert and R\"ollin~\cite{rr}. 
\end{proof}

\vspace{12pt}
Using normal comparison tools (i.e. results such as Slepian's Lemma, for comparing the behaviour of weakly correlated Gaussian random variables with independent Gaussians) one can readily establish the following assertion, that the maximum of $n$ standard Gaussians whose correlations are small is usually at least $\approx \sqrt{2\log n}$.

\begin{normcomp1}
Suppose that $n \geq 2$, and that $\epsilon > 0$ is sufficiently small (i.e. less than a certain small absolute constant). Let $X_1, ..., X_n$ be mean zero, variance one, jointly normal random variables, and suppose $\E X_i X_j \leq \epsilon$ whenever $i \neq j$.

Then for any $100\epsilon \leq \delta \leq 1/100$ (say), we have
$$ \p\left(\max_{1 \leq i \leq n} X_i \leq \sqrt{(2-\delta)\log n} \right) \ll e^{-\Theta(n^{\delta/20}/\sqrt{\log n})} + n^{-\delta^{2}/50\epsilon} . $$
\end{normcomp1}

\begin{proof}[Proof of Normal Comparison Result 1]
This is a special case of Normal Comparison Result 2 of Harper and Lamzouri~\cite{harperlamzouri}.
\end{proof}

\subsection{Number theory and harmonic analysis}
For our conditional variance and covariance calculations, we will need to introduce and analyse contour integrals involving Dirichlet polynomials (both random and deterministic). Our first tool for doing this, which is now widely used in problems involving random multiplicative functions, is the following version of Parseval's identity for Dirichlet series.

\begin{harman1}[See (5.26) in sec. 5.1 of Montgomery and Vaughan~\cite{mv}]
Let $(a_n)_{n=1}^{\infty}$ be any sequence of complex numbers, and let $A(s) := \sum_{n=1}^{\infty} \frac{a_n}{n^s}$ denote the corresponding Dirichlet series, and $\sigma_c$ denote its abscissa of convergence. Then for any $\sigma > \max\{0,\sigma_c \}$, we have
$$ \int_{0}^{\infty} \frac{|\sum_{n \leq x} a_n |^2}{x^{1 + 2\sigma}} dx = \frac{1}{2\pi} \int_{-\infty}^{\infty} \left|\frac{A(\sigma + it)}{\sigma + it}\right|^2 dt . $$
\end{harman1}

We will also need a mean value theorem for Dirichlet polynomials supported on fairly large primes and prime powers. The following version is well suited to our application (cruder versions would certainly also suffice, since this theorem is needed to control ``off diagonal'' terms in a double integral in a part of the argument that isn't too delicate).

\begin{numth1}
Uniformly for any complex numbers $(a_n)_{n=1}^{\infty}$ and any $T \geq 1$, and with $\Lambda(n)$ denoting the von Mangoldt function, we have
$$ \int_{-T}^{T} \Biggl|\sum_{T^{1.01} \leq n \leq x} \frac{a_n \Lambda(n)}{n^{1+it}} \Biggr|^{2} dt \ll \sum_{T^{1.01} \leq n \leq x} \frac{|a_n|^2 \Lambda(n)}{n} . $$
\end{numth1}

\begin{proof}[Proof of Number Theory Result 1]
This follows in a standard way combining a Gallagher-type upper bound for the integral (which ultimately comes from inserting a smooth weight and expanding out, see e.g. Lemme 3.1 of Tenenbaum~\cite{tenenbaummultmean} for a suitable general statement) with a sieve upper bound for primes in short intervals (this is where the condition $n \geq T^{1.01}$ in the sum is important).
\end{proof}

In various places we will require pointwise estimates for sums over primes $p$. We record some useful statements in the following result.

\begin{numth2}
For any $100 \leq x \leq y$ (say), any $|\sigma| \leq 1/\log y$, and any $|t| \leq 1/\log y$, we have $\sum_{x < p \leq y} \frac{1}{p^{1+2\sigma+2it}} = \log\log y - \log\log x + O(1)$.

Furthermore, for any $t \neq 0$ we have
$$ \left|\sum_{x < p \leq y} \frac{1}{p^{1+it}} \right| \leq \frac{3}{|t|\log x} + O((1+|t|) e^{-c\sqrt{\log x}}) . $$
\end{numth2}

\begin{proof}[Proof of Number Theory Result 2]
The first statement follows by writing $\frac{1}{p^{1+2\sigma+2it}} = \frac{1}{p} e^{-2(\sigma + it)\log p} = \frac{1}{p} (1 + O((|\sigma| + |t|)\log p))$ for $x < p \leq y$, and applying the standard Mertens estimate $\sum_{x < p \leq y} \frac{1}{p} = \log\log y - \log\log x + O(1)$ and Chebychev estimate $\sum_{x < p \leq y} \frac{\log p}{p} \ll \log y$ (see e.g. Theorem 2.7 of Montgomery and Vaughan~\cite{mv}).

The second statement can be deduced by applying partial summation to the Prime Number Theorem with classical error term.
\end{proof}

\subsection{Moments of random Euler products}
As in many other works involving random multiplicative functions, having passed to studying contour integrals and Dirichlet polynomials we will require estimates for moments of random Euler products. Unlike in previous works (see e.g. section 2 of Harper~\cite{harperrmfhigh}, or the closely related computations of shifted moments of $L$-functions by Chandee~\cite{chandee} and Soundararajan and Young~\cite{soundyoung}), in places we will need to deal with products of a very large (growing) number of {\em different} Euler factors, with good explicit dependence on the number of terms. The following results will provide sufficient uniformity.

\begin{epres1}
If $f$ is a Steinhaus random multiplicative function, then uniformly for any $k \in \N$, any real $\alpha_1, ..., \alpha_k$, any real $100(1 + (\sum_{j=1}^{k} |\alpha_j|)^2 ) \leq x \leq y$, and any real $\sigma \geq - 1/\log y$ and $t_1, ..., t_k$, we have
\begin{eqnarray}
&& \E \prod_{j=1}^{k} \prod_{x < p \leq y} \left|1 - \frac{f(p)}{p^{1/2+\sigma + it_j}}\right|^{-2\alpha_j} \nonumber \\
& = & \exp\{\sum_{x < p \leq y} \frac{\sum_{j=1}^{k} \alpha_j^2 + 2\sum_{1 \leq j < l \leq k} \alpha_j \alpha_l \cos((t_l - t_j)\log p)}{p^{1 + 2\sigma}} + O(\frac{M}{\sqrt{x} \log x}) \} , \nonumber
\end{eqnarray}
where $M = M(\alpha_1, ..., \alpha_k) := \max\{\sum_{j=1}^{k} |\alpha_j| , (\sum_{j=1}^{k} |\alpha_j| )^3\}$.

In particular, for any real $200 \leq x \leq y$, any real $\sigma \geq -1/\log y$, and any real $t_1 , t_2$ (possibly equal), we have
\begin{eqnarray}
&& \E \prod_{\substack{x < p \\ \leq y}} \left|1 - \frac{f(p)}{p^{1/2+\sigma + it_1}}\right|^{-1} \left|1 - \frac{f(p)}{p^{1/2+\sigma + it_2}}\right|^{-1} = \exp\{ \sum_{\substack{x < p \\ \leq y}} \frac{1 + \cos((t_1 - t_2)\log p)}{2p^{1 + 2\sigma}} + O(\frac{1}{\sqrt{x} \log x}) \} \nonumber \\
& \ll & \sqrt{\frac{\log y}{\log x} \left(1 + \min\left\{\frac{\log y}{\log x}, \frac{1}{|t_1 - t_2|\log x} + \frac{\log^{2}(2+|t_1 - t_2|)}{\log x} \right\} \right)} . \nonumber
\end{eqnarray}
\end{epres1}

We remark that the explicit final bound here is included for ease of reference later, but will not be sharp if $\sigma$ is large. 

\begin{proof}[Proof of Euler Product Result 1]
The proof is conceptually the same as for other results of this kind, developing each term $\prod_{j=1}^{k} \left|1 - \frac{f(p)}{p^{1/2+\sigma + it_j}}\right|^{-2\alpha_j}$ as a series expansion, carefully computing its expectation, and then multiplying together the contributions from different primes using independence. Since the uniformity of the result will be important later, we will briefly present the main steps.

Using the Taylor expansion of the logarithm, we may rewrite
\begin{eqnarray}
&& \prod_{j=1}^{k} \left|1 - \frac{f(p)}{p^{1/2+\sigma + it_j}}\right|^{-2\alpha_j} = \exp\{-2\sum_{j=1}^{k} \alpha_j \Re\log(1 - \frac{f(p)}{p^{1/2+\sigma + it_j}}) \} \nonumber \\
& = & \exp\{ \frac{2 \sum_{j=1}^{k} \alpha_j \Re f(p) p^{-it_j}}{p^{1/2 + \sigma}} + \frac{ \sum_{j=1}^{k} \alpha_j \Re f(p)^2 p^{-2it_j}}{p^{1+2\sigma}} + O(\frac{\sum_{j=1}^{k} |\alpha_j|}{p^{3/2 + 3\sigma}}) \} . \nonumber
\end{eqnarray}
And if $y \geq p > x \geq 100 (\sum_{j=1}^{k} |\alpha_j| )^2$, then each term in the exponential here has size at most $2(\sum_{j=1}^{k} |\alpha_j|)/p^{1/2 + \sigma} = 2(\sum_{j=1}^{k} |\alpha_j|) e^{-\sigma \log p}/ p^{1/2} \leq e/5$. Therefore we may apply the series expansion of the exponential function, finding the above is
$$ = 1 +  \frac{2 \sum_{j=1}^{k} \alpha_j \Re f(p) p^{-it_j}}{p^{1/2 + \sigma}} + \frac{ \sum_{j=1}^{k} \alpha_j \Re f(p)^2 p^{-2it_j}}{p^{1+2\sigma}} + \frac{2(\sum_{j=1}^{k} \alpha_j \Re f(p) p^{-it_j})^2}{p^{1 + 2\sigma}} + O(\frac{M}{p^{3/2 + 3\sigma}}) . $$

Using the fact that (by symmetry) we have $\E \Re f(p) p^{-it} = \E \Re f(p)^2 p^{-2it} = 0$, as well as $\E \Re f(p) \Re f(p) p^{-it} = \cos(t\log p)/2$, we find the expectation of the above expression is
\begin{eqnarray}
& = & 1 + \frac{\sum_{j=1}^{k} \alpha_j^2 + 2\sum_{1 \leq j < l \leq k} \alpha_j \alpha_l \cos((t_l - t_j)\log p)}{p^{1 + 2\sigma}} + O(\frac{M}{p^{3/2 + 3\sigma}}) \nonumber \\
& = & \exp\{\frac{\sum_{j=1}^{k} \alpha_j^2 + 2\sum_{1 \leq j < l \leq k} \alpha_j \alpha_l \cos((t_l - t_j)\log p)}{p^{1 + 2\sigma}} + O(\frac{M}{p^{3/2 + 3\sigma}}) \} . \nonumber
\end{eqnarray}
The first part of Euler Product Result 1 follows on using independence to multiply together the contributions from all primes $x < p \leq y$, noting that $\sum_{x < p \leq y} \frac{M}{p^{3/2 + 3\sigma}} \leq e^{3} M \sum_{x < p \leq y} \frac{1}{p^{3/2}} \ll \frac{M}{\sqrt{x} \log x}$ in the ``big Oh'' terms.

For the second part of Euler Product Result 1, the expression as an exponential of a prime number sum is just a special case of the first part of the Result, taking $k=2$ and $\alpha_1 = \alpha_2 =1/2$. Furthermore, if we temporarily set $\mathcal{X} := \max\{e^{1/|t_1 - t_2|}, e^{C\log^{2}(2+|t_1 - t_2|)}\}$, where $C$ is a suitable large constant, then the second part of Number Theory Result 2 (coupled with Abel summation to remove the factor $p^{-2\sigma}$) implies $\sum_{\mathcal{X} < p \leq y} \frac{\cos((t_1 - t_2)\log p)}{2p^{1+2\sigma}} \ll 1$. So we can upper bound our overall sum over primes $x < p \leq y$ by $\sum_{x < p \leq y} \frac{1 + \textbf{1}_{p \leq \mathcal{X}}}{2p^{1+2\sigma}} + O(1) \leq \sum_{x < p \leq y} \frac{1 + \textbf{1}_{p \leq \mathcal{X}}}{2p^{1-2/\log y}} + O(1)$, where $\textbf{1}$ denotes the indicator function. Applying the first part of Number Theory Result 2 to estimate this yields the final bound of Euler Product Result 1.
\end{proof}

\begin{epres2}
If $f$ is a Rademacher random multiplicative function, then uniformly for any $k \in \N$, any real $\alpha_1, ..., \alpha_k$, any real $100(1 + (\sum_{j=1}^{k} |\alpha_j|)^2 ) \leq x \leq y$, and any real $\sigma \geq - 1/\log y$ and $t_1, ..., t_k$, we have
\begin{eqnarray}
\E \prod_{j=1}^{k} \prod_{\substack{x < p \\ \leq y}} \left|1 + \frac{f(p)}{p^{1/2+\sigma + it_j}}\right|^{2\alpha_j} & = & \exp\{O(\frac{M}{\sqrt{x} \log x}) + \sum_{\substack{x < p \\ \leq y}} \Biggl( \frac{\sum_{j=1}^{k} \alpha_j^2 + \sum_{j=1}^{k} (\alpha_{j}^{2} - \alpha_j) \cos(2t_j \log p)}{p^{1 + 2\sigma}} + \nonumber \\
&& + \frac{2\sum_{1 \leq j < l \leq k} \alpha_j \alpha_l (\cos((t_l - t_j)\log p) + \cos((t_l + t_j)\log p))}{p^{1 + 2\sigma}} \Biggr) \} , \nonumber
\end{eqnarray}
where $M = M(\alpha_1, ..., \alpha_k) := \max\{\sum_{j=1}^{k} |\alpha_j| , (\sum_{j=1}^{k} |\alpha_j| )^3\}$.

In particular, for any real $200 \leq x \leq y$, any real $\sigma \geq -1/\log y$, and any real $t_1, t_2$ (possibly equal), we have
\begin{eqnarray}
&& \E \prod_{x < p \leq y} \left|1 + \frac{f(p)}{p^{1/2+\sigma + it_1}}\right| \left|1 + \frac{f(p)}{p^{1/2+\sigma + it_2}}\right| \nonumber \\
& \ll & \exp\{ \sum_{x < p \leq y} \frac{1 + \cos((t_1 - t_2)\log p) + \cos((t_1 + t_2)\log p) - \frac{1}{2}(\cos(2t_1 \log p) + \cos(2t_2\log p))}{2p^{1 + 2\sigma}} \} \nonumber \\
& \ll & \sqrt{\frac{\log y}{\log x} \left(1 + \min\left\{\frac{\log y}{\log x}, \frac{1}{|t_1 - t_2|\log x} + \frac{1}{|t_1 + t_2|\log x} + \frac{\log^{2}(2+|t_1| + |t_2|)}{\log x} \right\} \right)} . \nonumber
\end{eqnarray}
\end{epres2}

\begin{proof}[Proof of Euler Product Result 2]
The proof works in the same way as for Euler Product Result 1. Again, we outline the main steps.

Using the Taylor expansion of the logarithm, we find
\begin{eqnarray}
&& \prod_{j=1}^{k} \left|1 + \frac{f(p)}{p^{1/2+\sigma + it_j}}\right|^{2\alpha_j} = \exp\{2\sum_{j=1}^{k} \alpha_j \Re\log(1 + \frac{f(p)}{p^{1/2+\sigma + it_j}}) \} \nonumber \\
& = & \exp\{ \frac{2 \sum_{j=1}^{k} \alpha_j \Re f(p) p^{-it_j}}{p^{1/2 + \sigma}} - \frac{ \sum_{j=1}^{k} \alpha_j \Re f(p)^2 p^{-2it_j}}{p^{1+2\sigma}} + O(\frac{\sum_{j=1}^{k} |\alpha_j|}{p^{3/2 + 3\sigma}}) \} . \nonumber
\end{eqnarray}
As in Euler Product Result 1, using the series expansion of the exponential we can rewrite this as
$$ 1 +  \frac{2 \sum_{j=1}^{k} \alpha_j \Re f(p) p^{-it_j}}{p^{1/2 + \sigma}} - \frac{ \sum_{j=1}^{k} \alpha_j \Re f(p)^2 p^{-2it_j}}{p^{1+2\sigma}} + \frac{2(\sum_{j=1}^{k} \alpha_j \Re f(p) p^{-it_j})^2}{p^{1 + 2\sigma}} + O(\frac{M}{p^{3/2 + 3\sigma}}) . $$

In the Rademacher case we have $f(p)^2 \equiv 1$, but we still have $\E \Re f(p) p^{-it} = \E f(p) \Re p^{-it} = \cos(t\log p) \E f(p) = 0$, so a little calculation shows the expectation of the above expression is
\begin{eqnarray}
& = & 1 - \frac{\sum_{j=1}^{k} \alpha_j \cos(2t_j \log p)}{p^{1+2\sigma}} + \nonumber \\
&& + \frac{2(\sum_{j=1}^{k} \alpha_j^2 \cos^{2}(t_j \log p) + 2\sum_{1 \leq j < l \leq k} \alpha_j \alpha_l \cos(t_j \log p) \cos(t_l \log p))}{p^{1 + 2\sigma}} + O(\frac{M}{p^{3/2 + 3\sigma}}) . \nonumber
\end{eqnarray}
Using the trigonometric identity $\cos(t_j \log p) \cos(t_l \log p) = (1/2)(\cos((t_j - t_l) \log p) + \cos((t_j + t_l) \log p))$, and then multiplying together the contributions from different primes using independence, the first part of Euler Product Result 2 follows.

For the second part we proceed similarly as in Euler Product Result 1, first taking $k=2$ and $\alpha_1 = \alpha_2 =1/2$ in the preceding calculations. If we now temporarily set $\mathcal{X} := \max\{e^{1/|t_1 - t_2|}, e^{1/|t_1 + t_2|}, e^{C\log^{2}(2+|t_1| + |t_2|)}\}$, where $C$ is a suitable large constant, then the second part of Number Theory Result 2 (coupled with Abel summation to remove the factor $p^{-2\sigma}$) implies that $\sum_{\mathcal{X} < p \leq y} \frac{\cos((t_1 - t_2)\log p) + \cos((t_1 + t_2)\log p)}{2p^{1 + 2\sigma}} \ll 1$ and that $\sum_{\mathcal{X} < p \leq y} \frac{- \frac{1}{2}(\cos(2t_1 \log p) + \cos(2t_2\log p))}{2p^{1 + 2\sigma}}$ is {\em at most} a large constant. (Note that if $|t_1|$ or $|t_2|$ is small, in particular if $e^{1/|t_1|} > \mathcal{X}$ or $e^{1/|t_2|} > \mathcal{X}$, then the second sum might be large and {\em negative}, but it is always bounded above by a constant.) Meanwhile, using the identity $\frac{1}{2}(\cos(2t_1 \log p) + \cos(2t_2\log p)) = \cos((t_1 - t_2)\log p) \cos((t_1 + t_2)\log p)$ we see that the numerator of the whole prime number sum in the exponential may be rewritten as
$$ 2 - (1-\cos((t_1 + t_2)\log p))(1 - \cos((t_1 - t_2)\log p)) , $$
which is clearly $\leq 2$. So once again we can upper bound that sum by $\sum_{x < p \leq y} \frac{1 + \textbf{1}_{p \leq \mathcal{X}}}{2p^{1+2\sigma}} + O(1) \leq \sum_{x < p \leq y} \frac{1 + \textbf{1}_{p \leq \mathcal{X}}}{2p^{1-2/\log y}} + O(1)$, for which the first part of Number Theory Result 2 yields an acceptable bound.
\end{proof}

\subsection{Multiplicative chaos}\label{subsecmultchaos}
As remarked in the Introduction, it is now understood that the distribution of $\sum_{n \leq x} f(n)$ is closely connected, via {\em integrals} of random Euler products, to the probabilistic theory of (critical) multiplicative chaos. See the introduction to the author's paper~\cite{harperrmflow} for a detailed discussion of this. Here we develop some key points of the theory in a form suitable for our current purposes.

\vspace{12pt}
Let $X$ be large. For each $0 \leq k \leq \log\log X$ and complex $s$ with $\Re(s) > 0$, we shall write $F_{k}(s)$ for the partial Euler product corresponding to our random multiplicative function on $X^{e^{-k}}$-smooth numbers. More explicitly, if $f(n)$ is a Steinhaus random multiplicative function then $F_{k}(s) = \prod_{p \leq X^{e^{-k}}} (1 - \frac{f(p)}{p^s})^{-1}$, whereas if $f(n)$ is a Rademacher random multiplicative function (so only supported on squarefree numbers) then $F_{k}(s) = \prod_{p \leq X^{e^{-k}}} (1 + \frac{f(p)}{p^s})$. For concision we usually write $F(s)$ in place of $F_{0}(s)$.

To understand the typical or the distributional behaviour of integrals like $\int_{-1/2}^{1/2} |F(1/2+it)|^2 dt$, the most fundamental insight from critical multiplicative chaos is that one should restrict to the case where the partial products $F_{k}(1/2+it)$ (or closely related objects) satisfy certain size bounds at a net of points $t$, and moreover that one can choose these bounds so that it is highly probable they will be satisfied. More explicitly, for each $t \in \R$ we create a sequence of increasingly coarse approximations to $t$, by setting $t(-1) = t$ and then (for each $0 \leq j \leq \log\log X - 1$) setting
$$ t(j) := \max\{u \leq t(j-1) : u = \frac{n}{((\log X)/e^j) \log((\log X)/e^j)} \; \text{for some} \; n \in \Z \} . $$
These $t(j)$ will provide our net of points, and one can readily prove the following.

\begin{chaos1}
Let $f(n)$ be a Steinhaus or Rademacher random multiplicative function. Uniformly for all large $X$, all $N \in \R$, and all $W \geq 0$, we have
\begin{eqnarray}
&& \p\left(|F_{j}(\frac{1}{2} + it(j))| \leq \frac{\log X}{e^j} (\log\log X - j)^2 e^W \;\;\; \forall \; 0 \leq j \leq \log\log X - 1, \forall \; |t-N| \leq 1/2 \right) \nonumber \\
& \geq & 1 - O(e^{-2W}) . \nonumber
\end{eqnarray}
\end{chaos1}

\begin{proof}[Proof of Multiplicative Chaos Result 1]
This is essentially a more straightforward version of the proof of Key Proposition 2 of Harper~\cite{harperrmflow}. By the union bound, the probability that the given event fails to occur is at most
$$ \sum_{0 \leq j \leq \log\log X - 1} \sum_{t(j)} \p\left(|F_{j}(1/2 + it(j))| > \frac{\log X}{e^j} (\log\log X - j)^2 e^W \right) , $$
where the inner sum is over all possible values $t(j)$ that can arise for $|t-N| \leq 1/2$ (of which there are $\ll ((\log X)/e^j ) \log((\log X)/e^j)$). By Chebychev's inequality, this is all
$$ \leq \frac{1}{e^{2W}} \sum_{0 \leq j \leq \log\log X - 1} \sum_{t(j)} \frac{1}{((\log X)/e^j)^2 (\log\log X - j)^4} \E|F_{j}(1/2 + it(j))|^2 . $$
Finally, since we have $\E|F_{j}(1/2 + it(j))|^2 \ll \frac{\log X}{e^j}$ (by the second part of Euler Product Results 1 and 2) the result follows on performing the sums over $t(j)$ and over $j$.
\end{proof}

For our work in section \ref{secmainresult} here (and unlike in the paper~\cite{harperrmflow}), we need to upgrade the bound $\frac{\log X}{e^j} (\log\log X - j)^2 e^W$ permitted by Multiplicative Chaos Result 1 to something stronger, showing we can disregard values of $|F(1/2+it)|$ even where the partial products are a bit less than $\frac{\log X}{e^j} (\log\log X - j)^2$. To do this we shall use the following expectation estimate, which encodes the fact that the most typical way for all the partial products to remain below $\frac{\log X}{e^j} (\log\log X - j)^2$ (in fact for simplicity we work with a larger barrier $\frac{\log X}{e^j} (\log\log X)^5$) is for them to evolve on a ``path'' that is {\em significantly} below this. 
\begin{chaos2}
Let $f(n)$ be a Steinhaus or Rademacher random multiplicative function. Then uniformly for all large $X$ and all $\frac{1}{(\log\log X)^{1000}} \leq |t| \leq \log^{1000}X$, say (this restriction on $t$ is only needed in the Rademacher case), the following is true: if we let $\mathcal{D}(t)$ denote the event that
$$ |F_{j}(1/2 + it(j))| \leq \frac{\log X}{e^j} (\log\log X)^5 \;\;\; \forall \; 0 \leq j \leq \log\log X - 1 ; $$
and let $\mathcal{A}(t)$ denote the event that
$$ |F_{j}(\frac{1}{2} + it)| \leq \frac{\log X}{e^j (\log\log X)^{1000}} \;\; \forall \; j \leq 0.99\log\log X , \;\; \text{and} \;\; |F_{j}(\frac{1}{2} + it)| \leq \frac{\log X}{e^j} (\log\log X)^6 \;\; \forall \; j , $$
say; then (with $\textbf{1}$ denoting the indicator function) we have
$$ \E |F(1/2+it)|^2 \textbf{1}_{\mathcal{D}(t)} \textbf{1}_{\mathcal{A}(t) \; \text{fails}} \ll \log X \frac{(\log\log\log X)^7}{\log\log X} . $$
\end{chaos2}

The key feature of this bound is that it is rather smaller than $\frac{\log X}{\sqrt{\log\log X}}$, the rough size of the Euler product integrals we shall deal with later.

\begin{proof}[Proof of Multiplicative Chaos Result 2]
Firstly we note the upper bound
$$ \E |F(1/2+it)|^2 \textbf{1}_{\mathcal{D}(t)} \textbf{1}_{\mathcal{A}(t) \; \text{fails}} \leq \E |F(1/2+it)|^2 \textbf{1}_{\mathcal{D}(t)} \textbf{1}_{\mathcal{D}^{*}(t) \; \text{fails}} + \E |F(1/2+it)|^2 \textbf{1}_{\mathcal{D}^{*}(t)} \textbf{1}_{\mathcal{A}(t) \; \text{fails}} , $$
where $\mathcal{D}^{*}(t)$ is the event that
$$ \frac{1}{\log\log X} \leq |F_{j}(1/2 + it)| \leq \frac{\log X}{e^j} (\log\log X)^6 \;\;\; \forall \; 0 \leq j \leq \log\log X - 1 . $$
Now if $\mathcal{D}(t)$ holds but $\mathcal{D}^{*}(t)$ fails, then there must exist some $0 \leq j \leq \log\log X - 1$ for which $|F_{j}(1/2 + it(j))| \leq \frac{\log X}{e^j} (\log\log X)^5$ but either $|F_{j}(1/2 + it)| > \frac{\log X}{e^j} (\log\log X)^6$, or $|F_{j}(1/2 + it)| < 1/\log\log X$. Thus we have an upper bound
\begin{eqnarray}
\E |F(\frac{1}{2}+it)|^2 \textbf{1}_{\mathcal{D}(t)} \textbf{1}_{\mathcal{D}^{*}(t) \; \text{fails}} & \leq & \sum_{0 \leq j \leq \log\log X - 1} \frac{1}{(\log\log X)^2} \E |F(\frac{1}{2}+it)|^2 \frac{|F_{j}(1/2 + it)|^2}{|F_{j}(1/2 + it(j))|^2} \nonumber \\
&& + \sum_{0 \leq j \leq \log\log X - 1} \frac{1}{(\log\log X)^2} \E \frac{|F(1/2 + it)|^2}{|F_{j}(1/2 + it)|^2} , \nonumber
\end{eqnarray}
and the second sum here is easily seen to be $\ll \sum_{0 \leq j \leq \log\log X - 1} \frac{1}{(\log\log X)^2} e^j \ll \frac{\log X}{(\log\log X)^2}$ using Euler Product Results 1 and 2.

It will also be fairly straightforward to handle the first sum in the previous display. First writing the argument in the Steinhaus case, using the independence of the $f(p)$ for different primes $p$ we can factorise
$$ \E |F(\frac{1}{2}+it)|^2 \frac{|F_{j}(\frac{1}{2} + it)|^2}{|F_{j}(\frac{1}{2} + it(j))|^2} = \E \prod_{X^{e^{-j}} < p \leq X} |1 - \frac{f(p)}{p^{1/2+it}}|^{-2} \cdot \E \prod_{p \leq X^{e^{-j}}} \frac{|1 - \frac{f(p)}{p^{1/2+it}}|^{-4}}{|1 - \frac{f(p)}{p^{1/2+it(j)}}|^{-2}} , $$
and by Euler Product Result 1 this is $\ll \exp\{\sum_{X^{e^{-j}} < p \leq X} \frac{1}{p} + \sum_{p \leq X^{e^{-j}}} \frac{5 - 4\cos((t-t(j))\log p)}{p} \} = \exp\{\sum_{p \leq X} \frac{1}{p} + O(\sum_{p \leq X^{e^{-j}}} \frac{(t-t(j))^2 \log^{2}p}{p}) \}$. Here by construction we have $|t-t(j)| \ll \frac{1}{((\log X)/e^j) \log((\log X)/e^j)}$, so Chebychev's bounds imply the ``big Oh'' term is $\ll (t-t(j))^2 (\frac{\log X}{e^j})^2 \ll 1$, and the whole exponential is $\ll \log X$. Overall we deduce that $\E |F(1/2+it)|^2 \textbf{1}_{\mathcal{D}(t)} \textbf{1}_{\mathcal{D}^{*}(t) \; \text{fails}} \ll \frac{\log X}{\log\log X}$, which is an acceptable bound. In the Rademacher case, one obtains the same bound after applying Euler Product Result 2, noting that the contribution from all the terms $\frac{2\cos(2t\log p)}{p}$ and $\frac{2\cos(2t(j)\log p)}{p}$ there cancels (up to an overall $O(1)$ term) with the contribution $\frac{-4\cos((t+t(j))\log p)}{p}$.

Next, if $\mathcal{D}^{*}(t)$ holds but $\mathcal{A}(t)$ fails then there must exist some $J \leq 0.99\log\log X$ for which $|F_{J}(1/2 + it)| > \frac{\log X}{e^J (\log\log X)^{1000}}$. Thus we can upper bound $\E |F(1/2+it)|^2 \textbf{1}_{\mathcal{D}^{*}(t)} \textbf{1}_{\mathcal{A}(t) \; \text{fails}}$ by $\sum_{0 \leq J \leq 0.99\log\log X} \E |F(1/2+it)|^2 \textbf{1}_{\mathcal{A}_{J}(t)} \textbf{1}_{\mathcal{B}_{J}(t)}$, where $\mathcal{A}_{J}(t)$ is the event that
$$ \frac{1}{\log\log X} \leq |F_{j}(\frac{1}{2} + it)| \leq \frac{\log X}{e^j} (\log\log X)^6 \;\;\; \forall \; J \leq j , \;\;\; \text{and} \;\;\; |F_{J}(\frac{1}{2} + it)| > \frac{\log X}{e^J (\log\log X)^{1000}} , $$
and where $\mathcal{B}_{J}(t)$ is the event that (writing the Euler products for the Steinhaus case)
$$ \left(\frac{\log X}{e^J}\right)^{-1} \frac{1}{(\log\log X)^7} \leq \prod_{X^{e^{-J}} < p \leq X^{e^{-j}}} |1 - \frac{f(p)}{p^{1/2+it}}|^{-1} \leq e^{J-j} (\log\log X)^{1006} \;\;\; \forall \; 0 \leq j \leq J - 1 . $$
Note that the lower bound in $\mathcal{B}_{J}(t)$ is imposed by the fact that $|F_{j}(1/2 + it)| \geq 1/\log\log X$ for all $j$ and $|F_{J}(1/2 + it)| \leq \frac{\log X}{e^J} (\log\log X)^6$, whenever $\mathcal{D}^{*}(t)$ holds.

Since the events $\mathcal{A}_{J}(t)$ and $\mathcal{B}_{J}(t)$ involve the $f(p)$ for disjoint sets of primes $p$, they are independent and in the sum $\sum_{0 \leq J \leq 0.99\log\log X} \E |F(1/2+it)|^2 \textbf{1}_{\mathcal{A}_{J}(t)} \textbf{1}_{\mathcal{B}_{J}(t)}$ we can factorise
$$ \E |F(1/2+it)|^2 \textbf{1}_{\mathcal{A}_{J}(t)} \textbf{1}_{\mathcal{B}_{J}(t)} = \E |F_{J}(1/2+it)|^2 \textbf{1}_{\mathcal{A}_{J}(t)} \cdot \E \prod_{X^{e^{-J}} < p \leq X} |1 - \frac{f(p)}{p^{1/2+it}}|^{-2} \textbf{1}_{\mathcal{B}_{J}(t)} . $$
Now in the Steinhaus case, combining the second part of Lemma 4 of Harper~\cite{harperrmflow} (with the length $x^{1/e}$ of the Euler product there replaced by $X$, with $\sigma = 0$, and with $t$ shifted to zero using translation invariance in law for Steinhaus Euler products) with Probability Result 1 of Harper~\cite{harperrmflow} and then Euler Product Result 1 yields that
$$ \E \prod_{X^{e^{-J}} < p \leq X} |1 - \frac{f(p)}{p^{1/2+it}}|^{-2} \textbf{1}_{\mathcal{B}_{J}(t)} \ll \frac{\log\log\log X}{1 + \sqrt{J}} \E \prod_{X^{e^{-J}} < p \leq X} |1 - \frac{f(p)}{p^{1/2+it}}|^{-2} \ll \frac{\log\log\log X}{1 + \sqrt{J}} e^J . $$
(We remark that the lower bound condition in the definition of $\mathcal{B}_{J}(t)$ is only important here so that the size conditions in Lemma 4 of Harper~\cite{harperrmflow} are satisfied, and after that can be discarded since we are just seeking an upper bound for the expectation.)

To handle the event $\mathcal{A}_{J}(t)$, we need one final technical manipulation to remove the smallest primes so that Lemma 4 of Harper~\cite{harperrmflow} will again be applicable. Thus we let $\mathcal{J}_{0} := \lfloor \log\log X - (\log\log\log X)^2 \rfloor$, so that $e^{(\log\log\log X)^2} \leq \frac{\log X}{e^{\mathcal{J}_0}} \leq e^{(\log\log\log X)^2 + 1}$. Then we have $\E |F_{J}(1/2+it)|^2 \textbf{1}_{\mathcal{A}_{J}(t)} \leq \E |F_{J}(1/2+it)|^2 \textbf{1}_{\mathcal{C}_{J}(t)}$, where $\mathcal{C}_{J}(t)$ is the event that
$$ \frac{1}{e^{(\log\log\log X)^2 + 1} (\log\log X)^7} \leq \prod_{X^{e^{-\mathcal{J}_{0}}} < p \leq X^{e^{-j}}} |1 - \frac{f(p)}{p^{1/2+it}}|^{-1} \leq \frac{\log X}{e^j} (\log\log X)^7 \;\;\; \forall \; J \leq j \leq \mathcal{J}_{0} - 1 , $$
{\em and} $\prod_{X^{e^{-\mathcal{J}_{0}}} < p \leq X^{e^{-J}}} |1 - \frac{f(p)}{p^{1/2+it}}|^{-1} > \frac{\log X}{e^{J+(\log\log\log X)^2 + 1} (\log\log X)^{1006}}$. (These bounds are obtained by comparing the $\mathcal{A}_{J}(t)$ restrictions on $|F_{J}(1/2+it)|$, on $|F_{\mathcal{J}_0}(1/2+it)|$, and on $|F_{j}(1/2+it)|$ for general $J \leq j \leq \mathcal{J}_{0} - 1$.) Now Lemma 4 of Harper~\cite{harperrmflow} may be used, since the upper and lower size bounds in $\mathcal{C}_{J}(t)$ satisfy the relevant conditions in terms of $\frac{\log X}{e^j}$, and it implies that the ratio $\frac{\E |F_{J}(1/2+it)|^2 \textbf{1}_{\mathcal{C}_{J}(t)}}{\E |F_{J}(1/2+it)|^2}$ is at most a constant times
$$ \p( \max_{k \leq \mathcal{J}_{0} - J} \sum_{m=1}^{k} G_m \leq (\log\log\log X)^2 + O(\log\log\log X) , \; \text{and} \; \sum_{m=1}^{\mathcal{J}_{0} - J} G_m \geq -1006\log\log\log X - O(1) ) , $$
where the $G_m$ are independent Gaussian random variables with mean zero and variance $\sum_{X^{e^{-\mathcal{J}_{0} + m - 1}} < p \leq X^{e^{-\mathcal{J}_{0} + m}}} \frac{1}{2p} = 1/2 + o(1)$. Using a suitable form of the Ballot Theorem on Gaussian random walks, e.g. Probability Result 1 of Harper~\cite{harperpartition}, one finds this probability is $\ll \frac{(\log\log\log X)^2}{\sqrt{\mathcal{J}_{0} - J}} (\frac{(\log\log\log X)^2}{\sqrt{\mathcal{J}_{0} - J}})^2 \ll \frac{(\log\log\log X)^6}{(\log\log X)^{3/2}}$, and so overall
$$ \E |F_{J}(1/2+it)|^2 \textbf{1}_{\mathcal{C}_{J}(t)} \ll \frac{(\log\log\log X)^6}{(\log\log X)^{3/2}} \E |F_{J}(1/2+it)|^2 \ll \frac{(\log\log\log X)^6}{(\log\log X)^{3/2}} \frac{\log X}{e^J} . $$

Putting things together, we have shown that $\E |F(1/2+it)|^2 \textbf{1}_{\mathcal{D}^{*}(t)} \textbf{1}_{\mathcal{A}(t) \; \text{fails}}$ is
$$ \ll \sum_{0 \leq J \leq 0.99\log\log X} \frac{(\log\log\log X)^6}{(\log\log X)^{3/2}} \frac{\log X}{e^J} \frac{\log\log\log X}{1 + \sqrt{J}} e^J \ll \frac{(\log\log\log X)^7}{\log\log X} \log X , $$
which is acceptable. 

In the Rademacher case, one proceeds in the same way to analyse the Rademacher Euler products (and with the events $\mathcal{A}_{J}(t), \mathcal{B}_{J}(t), \mathcal{C}_{J}(t)$ now defined in terms of the Rademacher products). The term $\E \prod_{X^{e^{-J}} < p \leq X} |1 + \frac{f(p)}{p^{1/2+it}}|^{2} \textbf{1}_{\mathcal{B}_{J}(t)}$ may be bounded using Lemma 5 and Probability Result 1 of Harper~\cite{harperrmflow}, yielding the same bound $\frac{\log\log\log X}{1 + \sqrt{J}} e^J$ as in the Steinhaus case. Our condition that $\frac{1}{(\log\log X)^{1000}} \leq |t| \leq \log^{1000}X$ is used here (although at this point much weaker conditions would suffice), along with the observation that $X^{e^{-J}} \geq e^{\log^{0.01}X}$, to satisfy the assumptions of Lemma 5. Similarly, the term $\E |F_{J}(1/2+it)|^{2} \textbf{1}_{\mathcal{A}_{J}(t)}$ may be upper bounded by $\E |F_{J}(1/2+it)|^{2} \textbf{1}_{\mathcal{C}_{J}(t)}$, and then bounded using Lemma 5 of Harper~\cite{harperrmflow} and the Ballot Theorem. Again, one needs the condition $\frac{1}{(\log\log X)^{1000}} \leq |t| \leq \log^{1000}X$ here, along with the observation that $X^{e^{-\mathcal{J}_0}} \geq  e^{e^{(\log\log\log X)^2}}$, to satisfy the assumptions of Lemma 5. One obtains the same bound $\frac{(\log\log\log X)^6}{(\log\log X)^{3/2}} \frac{\log X}{e^J}$ as in the Steinhaus case, and so the same overall result.
\end{proof}

\vspace{12pt}
By calculating quantities like $\E \textbf{1}_{\mathcal{G}} \int_{-1/2}^{1/2} |F(1/2+it)|^2 dt$ for events $\mathcal{G}$ of the form discussed above (that restrict the size of partial products), and then applying H\"older's inequality to interpolate, one can prove bounds like the following.
\begin{chaos3}
Let $f(n)$ be a Steinhaus or Rademacher random multiplicative function. Then uniformly for all large $X$, all $0 \leq q \leq 1$, all $-1/\log X \leq \sigma \leq 1/\log^{0.01}X$ (say), and all $|N| \leq \log^{1000}X$ (say), we have
$$ \E\Biggl(\int_{N-1/2}^{N+1/2} |F(1/2 + \sigma + it)|^2 dt \Biggr)^q \ll (\log\log(|N| + 10))^q \Biggl(\frac{\min\{\log X, 1/|\sigma|\}}{1 + (1-q)\sqrt{\log\log X}} \Biggr)^q . $$

In the Steinhaus case, the factor $(\log\log(|N| + 10))^q$ may be omitted.
\end{chaos3}

\begin{proof}[Proof of Multiplicative Chaos Result 3]
We shall make some fairly simple reductions, and then the claimed bound will be read quite directly out of the author's paper~\cite{harperrmflow}. 

In the Steinhaus case, the distribution of the sequence $(f(n)n^{-it})_{n \in \N}$ is the same for any fixed $t \in \R$, so the distribution of $F(1/2 + \sigma + it)$ is invariant under shifts in $t$ and we only need to handle the case $N = 0$. Furthermore, if we temporarily let $\E^{(\sigma)}$ denote expectation conditional on the values $(f(p))_{p \leq \min\{X, e^{1/|\sigma|}\}}$, and let $F^{(\sigma)}(s) := \prod_{p \leq \min\{X, e^{1/|\sigma|}\}} (1 - \frac{f(p)}{p^s})^{-1}$, then using (the conditional form of) H\"older's inequality we find $\E(\int_{-1/2}^{1/2} |F(1/2 + \sigma + it)|^2 dt )^q = \E \E^{(\sigma)} (\int_{-1/2}^{1/2} |F(\frac{1}{2} + \sigma + it)|^2 dt)^q$ is
$$ \leq \E \Biggl(\int_{-1/2}^{1/2} \E^{(\sigma)} |F(\frac{1}{2} + \sigma + it)|^2 \Biggr)^q = \E \Biggl(\int_{-1/2}^{1/2} |F^{(\sigma)}(\frac{1}{2} + \sigma + it)|^2 \E \prod_{\min\{X, e^{1/|\sigma|}\} < p \leq X} |1 - \frac{f(p)}{p^{1/2+\sigma+it}}|^{-2} \Biggr)^q . $$
Using Euler Product Result 1, the expectation of the second product here has size $\ll \exp\{\sum_{\min\{X, e^{1/|\sigma|}\} < p \leq X} \frac{1}{p^{1+2\sigma}}\}$. And under the assumptions of Multiplicative Chaos Result 3, this sum over primes is empty unless $1/\log X \leq \sigma \leq 1/\log^{0.01}X$, in which case the sum is $\sum_{e^{1/\sigma} < p \leq X} \frac{1}{p^{1+2\sigma}} \leq \sigma \sum_{e^{1/\sigma} < p \leq X} \frac{\log p}{p^{1+2\sigma}} \ll \sigma \int_{e^{1/\sigma}}^{X} \frac{1}{t^{1+2\sigma}} dt \ll 1$ by Chebychev's estimates for primes.

It will now suffice to establish the claimed bound for $\E (\int_{-1/2}^{1/2} |F^{(\sigma)}(1/2 + \sigma + it)|^2 dt )^q$. In fact it will suffice to do this for $2/3 \leq q \leq 1$, since the result for smaller $q$ then follows by using H\"older's inequality to compare with the case $q=2/3$. This bound follows as in the proof of the upper bound part of Theorem 1 of Harper~\cite{harperrmflow}, taking $k=0$ there and with the length $x^{1/e}$ of the Euler product replaced by $\min\{X, e^{1/|\sigma|}\}$ (note that the shift $\sigma$ now makes no difference to any of the calculations, since the logarithm of the length of the Euler product $F^{(\sigma)}(s)$ is $\leq 1/|\sigma|$).

In the Rademacher case we cannot restrict to $N=0$, but can perform the same manipulations as before to replace $F(s)$ by $F^{(\sigma)}(s)$. The desired bound for $\E (\int_{N-1/2}^{N+1/2} |F^{(\sigma)}(1/2 + \sigma + it)|^2 dt )^q$ again follows as in the proof of the upper bound in Theorem 2 of Harper~\cite{harperrmflow}, noting that the restriction to $|N| \leq (\log\log x)^2$ imposed there may be relaxed to $|N| \leq \log^{1000}X$ (and in fact relaxed much more) without changing the proofs.
\end{proof}

The key feature here is the denominator $1 + (1-q)\sqrt{\log\log X}$. A direct application of H\"older's inequality, upper bounding the $q$-th moment of the integral by the $q$-th power of the first moment, would produce an analogous bound but without this saving.

\vspace{12pt}
We shall also require lower bounds for quantities like $\int_{-1/2}^{1/2} |F(1/2+it)|^2 dt$. The upper bound in Multiplicative Chaos Result 3 is sharp (see the author's paper~\cite{harperrmflow}), so comparing the order of the moments with $q=1/2$ and $q=3/4$ (say) implies that $\int_{-1/2}^{1/2} |F(1/2 + \sigma + it)|^2 dt \gg \frac{\min\{\log X, 1/|\sigma|\}}{\sqrt{\log\log X}}$ with positive probability. But we need a lower bound that holds with probability very close to 1, and deducing that will require more work. To avoid technical irritations that would arise in the Rademacher case when handling $t$ close to zero (because in that case the distribution of $(f(n)n^{-it})_{n \in \N}$ is not the same for all $t \in \R$), we shall actually focus on lower bounding $\int_{1/3}^{1/2} |F(1/2 + \sigma + it)|^2 dt$.

\begin{chaos4}
Let $f(n)$ be a Steinhaus or Rademacher random multiplicative function. Then uniformly for all large $X$, all $-1/\log X \leq \sigma \leq 1/\log^{0.01}X$ (say), and all $1 \leq W \leq (\log\log X)^{1/100}$ (say), we have
$$ \p\left(\int_{1/3}^{1/2} |F(1/2 + \sigma + it)|^2 dt \geq \frac{1}{e^{2.1 W}} \frac{\min\{\log X, 1/|\sigma|\}}{\sqrt{\log\log X}} \right) \geq 1 - O(e^{-0.1W}) . $$
\end{chaos4}

\begin{proof}[Proof of Multiplicative Chaos Result 4]
To obtain a probability bound that is close to 1, the key issue is to first extract the contribution to $F(1/2 + \sigma + it)$ from ``small'' primes, which doesn't vary much as $t$ varies so shouldn't be analysed in the same way as the large prime contribution (which ``averages out'' more over the integral). Simultaneously, we will extract the contribution from any primes larger than $e^{1/|\sigma|}$, which is insignificant (killed off by the shift by $\sigma$ in the real part) as in Multiplicative Chaos Result 3. Following these reductions, we will be able to analyse the contribution from the remaining ``medium'' primes by adapting the proof of Key Proposition 5 of Harper~\cite{harperrmflow}.  

More precisely, let $B$ denote a suitable large constant, and let us temporarily set $\mathcal{X} := \min\{X, e^{1/|\sigma|}\}$. Then we let $\mathcal{N}$ denote the largest natural number for which $\mathcal{X}^{e^{-\mathcal{N}}} \geq e^{e^{W+B}}$, noting that $\log(\mathcal{X}^{e^{-\mathcal{N}}}) \asymp e^W$ and that $\mathcal{N} = \log\log\mathcal{X} - W - B + O(1) \asymp \log\log X$. Finally, in the Steinhaus case we define $F^{(\text{med})}(s) := \prod_{\mathcal{X}^{e^{-\mathcal{N}}} < p \leq \mathcal{X}} (1 - \frac{f(p)}{p^s})^{-1}$, and define $\mathcal{L}$ to be the random subset of $[1/3,1/2]$ consisting of points $t$ for which
$$ \left(\frac{\log\mathcal{X}}{e^{j+W}} \right)^{-B} \leq \prod_{\mathcal{X}^{e^{-\mathcal{N}}} < p \leq \mathcal{X}^{e^{-j}}} |1 - \frac{f(p)}{p^{1/2+\sigma+it}}|^{-1} \leq \frac{\log\mathcal{X}}{e^{j+W}} (\log\log\mathcal{X} - j - W)^{-2} \;\;\; \forall \; 0 \leq j \leq \mathcal{N} - 1 . $$
Then by the Cauchy--Schwarz inequality, $\int_{1/3}^{1/2} |F(1/2 + \sigma + it)|^2 dt$ is always
\begin{equation}\label{chaoscsequation}
\geq \int_{\mathcal{L}} |F(1/2 + \sigma + it)|^2 dt \geq \frac{(\int_{\mathcal{L}} |F^{(\text{med})}(1/2 + \sigma + it)|^2 dt )^{2}}{\int_{\mathcal{L}} |F^{(\text{med})}(1/2 + \sigma + it)|^2 \cdot \prod_{\substack{p \leq \mathcal{X}^{e^{-\mathcal{N}}}, \\ \text{or} \; \mathcal{X} < p \leq X}} |1 - \frac{f(p)}{p^{1/2+\sigma+it}}|^{2} dt} .
\end{equation}
(Notice the exponent of the extra Euler factors in the denominator is 2, rather than -2.)

The definition of the random subset $\mathcal{L}$ only involves the $f(p)$ for primes occurring in $F^{(\text{med})}(s)$. Thus by independence of the different $f(p)$, the expectation of the denominator in \eqref{chaoscsequation} is $= \int_{1/3}^{1/2} \E \textbf{1}_{t \in \mathcal{L}} |F^{(\text{med})}(1/2 + \sigma + it)|^2 \cdot \E \prod_{\substack{p \leq \mathcal{X}^{e^{-\mathcal{N}}}, \\ \text{or} \; \mathcal{X} < p \leq X}} |1 - \frac{f(p)}{p^{1/2+\sigma+it}}|^{2} dt$, which is
\begin{eqnarray}
& \ll & \int_{1/3}^{1/2} \E \textbf{1}_{t \in \mathcal{L}} |F^{(\text{med})}(1/2 + \sigma + it)|^2 \cdot \exp\{\sum_{\substack{p \leq \mathcal{X}^{e^{-\mathcal{N}}}, \\ \text{or} \; \mathcal{X} < p \leq X}} \frac{1}{p^{1+2\sigma}}\} dt \nonumber \\
& \ll & \log(\mathcal{X}^{e^{-\mathcal{N}}}) \int_{1/3}^{1/2} \E \textbf{1}_{t \in \mathcal{L}} |F^{(\text{med})}(\frac{1}{2} + \sigma + it)|^2 dt \ll e^{W} \int_{1/3}^{1/2} \E \textbf{1}_{t \in \mathcal{L}} |F^{(\text{med})}(\frac{1}{2} + \sigma + it)|^2 dt . \nonumber
\end{eqnarray}
So by Markov's inequality, that denominator will be less than $e^{1.1 W} \int_{1/3}^{1/2} \E \textbf{1}_{t \in \mathcal{L}} |F^{(\text{med})}(1/2 + \sigma + it)|^2 dt$ with probability $\geq 1 - O(e^{-0.1W})$. We also note, for use in a short while, that
$$ \int_{1/3}^{1/2} \E \textbf{1}_{t \in \mathcal{L}} |F^{(\text{med})}(\frac{1}{2} + \sigma + it)|^2 dt \asymp \frac{1}{\sqrt{\mathcal{N}}} \exp\{\sum_{\mathcal{X}^{e^{-\mathcal{N}}} < p \leq \mathcal{X}} \frac{1}{p^{1+2\sigma}}\} \asymp \frac{\log\mathcal{X}}{e^{W} \sqrt{\log\log X}} , $$
which follows by combining Lemma 4 of Harper~\cite{harperrmflow} (with the length $x^{1/e}$ of the Euler product there replaced by $\mathcal{X}$, and with $t$ shifted to zero using translation invariance in law for Steinhaus Euler products) with Probability Result 2 of Harper~\cite{harperrmflow}.

To obtain a lower bound for the numerator in \eqref{chaoscsequation}, we shall compute the variance
\begin{eqnarray}
&& \E \left(\int_{\mathcal{L}} |F^{(\text{med})}(1/2 + \sigma + it)|^2 dt - \E \int_{\mathcal{L}} |F^{(\text{med})}(1/2 + \sigma + it)|^2 dt \right)^{2} \nonumber \\
& = & \E \left(\int_{\mathcal{L}} |F^{(\text{med})}(1/2 + \sigma + it)|^2 dt \right)^2 - \left(\E \int_{\mathcal{L}} |F^{(\text{med})}(1/2 + \sigma + it)|^2 dt \right)^{2} \nonumber \\
& = & \int \int \E \textbf{1}_{t \in \mathcal{L}} |F^{(\text{med})}(\frac{1}{2} + \sigma + it)|^2 \textbf{1}_{u \in \mathcal{L}} |F^{(\text{med})}(\frac{1}{2} + \sigma + iu)|^2 - \left(\E \int_{\mathcal{L}} |F^{(\text{med})}(\frac{1}{2} + \sigma + it)|^2 \right)^{2} . \nonumber
\end{eqnarray}
When $|t-u| > e^{-W/4}$, the second part of Lemma 7 of Harper~\cite{harperrmflow} (applied with $x^{1/e}$ replaced by $\mathcal{X}$, and $t$ replaced by $t-u$, and combined with the slicing argument from Lemma 4 of Harper~\cite{harperrmflow}) implies\footnote{Strictly speaking, to obtain the exact decorrelation estimate claimed here one should slightly change the definition of the random set $\mathcal{L}$, so that the event $t \in \mathcal{L}$ can be exactly sliced up into events of the form treated in Lemma 7 without any error at the boundaries. Making this boundary change wouldn't affect any of our other calculations, but would make the definition of $\mathcal{L}$ much messier to write, so we leave this to the careful reader.} that $\E \textbf{1}_{t \in \mathcal{L}} |F^{(\text{med})}(\frac{1}{2} + \sigma + it)|^2 \textbf{1}_{u \in \mathcal{L}} |F^{(\text{med})}(\frac{1}{2} + \sigma + iu)|^2$ is
$$ = (1 + O(e^{-W/4})) \E \textbf{1}_{t \in \mathcal{L}} |F^{(\text{med})}(\frac{1}{2} + \sigma + it)|^2 \cdot \E \textbf{1}_{u \in \mathcal{L}} |F^{(\text{med})}(\frac{1}{2} + \sigma + iu)|^2 . $$
(This decorrelation estimate is where the fact that $F^{(\text{med})}(s)$ only involves primes larger than $e^{e^{W+B}}$, i.e. significantly larger than $e^{1/|u-t|}$, is crucial.) Meanwhile, working through the proof of Key Proposition 5 of Harper~\cite{harperrmflow} (with the length $x$ of the Euler product there replaced by $\mathcal{X}$, and with $V=1$ and $q=2/3$ say), systematically extracting the part of $|F^{(\text{med})}(\frac{1}{2} + \sigma + iu)|^2$ involving primes $\leq e^{1/|u-t|}$ and bounding this using the definition of the random subset $\mathcal{L}$, we find the contribution to the double integral from points $t,u$ with $|t-u| \leq e^{-W/4}$ is
\begin{eqnarray}
& \ll & \int_{1/3}^{1/2} \Biggl(\int_{|u-t| \leq 1/\log^{1/3}\mathcal{X}} \frac{\log^{2}\mathcal{X}}{e^{3W} (\log\log X)^4} \min\{\log\mathcal{X}, \frac{1}{|u-t|}\} du + \nonumber \\
&& + \int_{\substack{\frac{1}{\log^{1/3}\mathcal{X}} \\ < |u-t| \leq e^{-W}}} \frac{\log^{2}\mathcal{X} \cdot \log^{2}(\frac{2}{|u-t|})}{e^{3W} \log\log X \cdot |u-t| \log^{4}(\frac{2}{e^{W}|u-t|})} du + \int_{\substack{e^{-W} \\ < |u-t| \leq e^{-W/4}}} \frac{\log^{2}\mathcal{X} \cdot \log^{2}(\frac{2}{|u-t|})}{e^{2W} \log\log X} du \Biggr) dt \nonumber \\
& \ll & \frac{\log^{2}\mathcal{X}}{e^{3W} (\log\log X)^3} + \frac{\log^{2}\mathcal{X} \cdot W^2}{e^{3W} \log\log X} + \frac{\log^{2}\mathcal{X}}{e^{2W} \log\log X} \frac{W^2}{e^{W/4}} \ll \frac{\log^{2}\mathcal{X}}{e^{2W} \log\log X} \frac{W^2}{e^{W/4}} . \nonumber
\end{eqnarray}
So combining all the above calculations with the fact that $\int_{1/3}^{1/2} \E \textbf{1}_{t \in \mathcal{L}} |F^{(\text{med})}(\frac{1}{2} + \sigma + it)|^2 \asymp \frac{\log\mathcal{X}}{e^{W} \sqrt{\log\log X}}$, we have computed that the variance of $\int_{\mathcal{L}} |F^{(\text{med})}(1/2 + \sigma + it)|^2 dt$ is
$$ \ll e^{-W/4} \left(\E \int_{\mathcal{L}} |F^{(\text{med})}(\frac{1}{2} + \sigma + it)|^2 \right)^{2} + \frac{\log^{2}\mathcal{X}}{e^{2W} \log\log X} \frac{W^2}{e^{W/4}} \ll \frac{W^2}{e^{W/4}} \left(\E \int_{\mathcal{L}} |F^{(\text{med})}(\frac{1}{2} + \sigma + it)|^2 \right)^{2} . $$
Combining this good variance bound with Chebychev's inequality, we deduce that $\int_{\mathcal{L}} |F^{(\text{med})}(1/2 + \sigma + it)|^2 dt$ will be $\geq (1/2) \E \int_{\mathcal{L}} |F^{(\text{med})}(1/2 + \sigma + it)|^2 dt \gg \frac{\log\mathcal{X}}{e^W \sqrt{\log\log X}}$ with probability $\geq 1 - O(W^2 e^{-W/4})$.

Inserting our bounds for the numerator and denominator into \eqref{chaoscsequation}, we conclude as desired that with probability $\geq 1 - O(e^{-0.1W}) - O(W^2 e^{-W/4}) = 1 - O(e^{-0.1W})$, we have
$$ \int_{1/3}^{1/2} |F(1/2 + \sigma + it)|^2 dt \gg \frac{(\frac{\log\mathcal{X}}{e^W \sqrt{\log\log X}} )^{2}}{e^{1.1 W} \frac{\log\mathcal{X}}{e^W \sqrt{\log\log X}}} = \frac{\log\mathcal{X}}{e^{2.1W} \sqrt{\log\log X}} . $$

In the Rademacher case, one performs the analogous calculations with the Rademacher Euler products $\prod_{p} (1 + \frac{f(p)}{p^s})$. When computing the expectation of the denominator in the analogue of \eqref{chaoscsequation}, one obtains a factor $\exp\{\sum_{\substack{p \leq \mathcal{X}^{e^{-\mathcal{N}}}, \\ \text{or} \; \mathcal{X} < p \leq X}} \frac{1 + 2\cos(2t\log p)}{p^{1+2\sigma}}\}$ from Euler Product Result 2 in place of $\exp\{\sum_{\substack{p \leq \mathcal{X}^{e^{-\mathcal{N}}}, \\ \text{or} \; \mathcal{X} < p \leq X}} \frac{1}{p^{1+2\sigma}}\}$. However, when $1/3 \leq t \leq 1/2$ the second part of Number Theory Result 2 (along with Abel summation to remove the factor $p^{-2\sigma}$) implies the sum of all the terms $\frac{\cos(2t\log p)}{p^{1+2\sigma}}$ is $O(1)$, which doesn't alter anything. The other calculations also go through (using e.g. Lemma 5 of Harper~\cite{harperrmflow} instead of Lemma 4), up to changes of the same kind.

\end{proof}

% SECTION 3 %%%%%%%%%%%%%%%%%%%%%%%%%%%%%
\section{Proof of the main result}\label{secmainresult}
We are now ready to prove our main result, Theorem \ref{thmlocalas}. The proof is in four parts: firstly we shall use our Gaussian approximation and comparison tools, together with some of the multiplicative chaos calculations already performed in section \ref{sectools}, to complete the proof except for obtaining suitable bounds for the covariances of our random variables; then we shall use our analytic number theory tools, and a little calculation with random Euler products, to find a more convenient expression for those covariances; thirdly we shall use more of our multiplicative chaos calculations to show that (with high probability) we can impose a strong ``barrier'' condition inside the covariance expression; and finally we obtain acceptable covariance bounds by combining the barrier condition, moment estimates for random Euler products, and some new arguments.

\subsection{Reducing the proof to covariance estimates}\label{subsecmostofproof}
Firstly we work in the Steinhaus case. We are going to investigate the distribution of the sums $\sum_{n \leq x} f(n)$, for numbers $x$ of the form $X^{8/7} e^{2\pi r}$ and $0 \leq r \leq 2(\log X)/(21\pi)$. Thus all of these $x$-values satisfy $X^{8/7} \leq x \leq X^{8/7 + 4/21} = X^{4/3}$, and using multiplicativity we can break up the sums as
$$ \sum_{n \leq x} f(n) = \sum_{X < p \leq x} f(p) \sum_{m \leq x/p} f(m) + \sum_{\substack{n \leq x, \\ n \; \text{is} \; X \; \text{smooth}}} f(n) . $$

We do not wish to work much with the part of the sum over $X$-smooth numbers (i.e. numbers with all their prime factors $\leq X$). To remove this, we simply note that
\begin{eqnarray}
&& \E\#\biggl\{0 \leq r \leq \frac{2\log X}{21\pi} : |\sum_{\substack{n \leq X^{8/7} e^{2\pi r}, \\ n \; \text{is} \; X \; \text{smooth}}} f(n)| \geq \sqrt{X^{8/7} e^{2\pi r}} (\log\log X)^{0.01}\biggr\} \nonumber \\
& \leq & \sum_{0 \leq r \leq 2(\log X)/(21\pi)} \E \frac{|\sum_{\substack{n \leq X^{8/7} e^{2\pi r}, \\ n \; \text{is} \; X \; \text{smooth}}} f(n)|^2}{X^{8/7} e^{2\pi r} (\log\log X)^{0.02}} \ll \frac{\log X}{(\log\log X)^{0.02}} , \nonumber
\end{eqnarray}
using linearity of expectation and an easy second moment bound. So by Markov's inequality, with probability $\geq 1 - O((\log\log X)^{-0.02})$ there will exist a random set $\mathcal{R} \subseteq \{0 \leq r \leq 2(\log X)/(21\pi)\}$, with cardinality at least $2(\log X)/(22\pi) = (\log X)/(11\pi)$ (say), such that
$$ |\sum_{\substack{n \leq X^{8/7} e^{2\pi r}, \\ n \; \text{is} \; X \; \text{smooth}}} f(n)| < \sqrt{X^{8/7} e^{2\pi r}} (\log\log X)^{0.01} \;\;\;\;\;\ \forall \; r \in \mathcal{R} . $$
Note carefully that the sums over $X$-smooth numbers only depend on the random variables $(f(p))_{p \leq X}$. Thus the probability $1 - O((\log\log X)^{-0.02})$ is over all realisations of $(f(p))_{p \leq X}$, and the random set $\mathcal{R}$ only depends on $(f(p))_{p \leq X}$ (and in particular is independent of $(f(p))_{X < p \leq X^{4/3}}$).

\vspace{12pt}
Now we condition on the values $(f(p))_{p \leq X}$, and for $x = X^{8/7} e^{2\pi r}$ and $r \in \mathcal{R}$ we consider the distribution of $\sum_{X < p \leq x} f(p) \sum_{m \leq x/p} f(m)$. To simplify the analysis, we shall only investigate the real parts of these sums (if we can find a sum with real part significantly larger than $\sqrt{x} (\log\log X)^{0.01}$, then $\sum_{n \leq x} f(n)$ will have large real part and so large absolute value). Note also that $x/p \leq X^{4/3}/X = X^{1/3}$ in the inner sums, so {\em having conditioned on $(f(p))_{p \leq X}$ these inner sums become fixed}.

To simplify the writing, we let $\tilde{\p}$ denote probabilities conditional on the values $(f(p))_{p \leq X}$, and let $\tilde{\E}$ denote the corresponding conditional expectation, and let $\mathcal{X} := \{X^{8/7} e^{2\pi r} : r \in \mathcal{R}\}$ (so $\mathcal{X}$ is a random set, but only depends on $(f(p))_{p \leq X}$ and becomes fixed having first conditioned on those values). Then applying Normal Approximation Result 1 to the sums $\Re \frac{1}{\sqrt{x}} \sum_{X < p \leq x} f(p) \sum_{m \leq x/p} f(m)$ under the measure $\tilde{\p}$, we find that for any $u \in \R$ and any small $\eta > 0$, the conditional probability $\tilde{\p}(\max_{x \in \mathcal{X}} \Re \frac{1}{\sqrt{x}} \sum_{X < p \leq x} f(p) \sum_{m \leq x/p} f(m) \leq u )$ is
\begin{eqnarray}\label{gaussapproxdisplay}
& \leq & \p(\max_{x \in \mathcal{X}} Z_x \leq u + \eta ) + O\Biggl(\frac{1}{\eta^{2}} \sum_{x,y \in \mathcal{X}} \sqrt{\sum_{X < p \leq X^{4/3}} \frac{1}{xy}|\sum_{m \leq x/p} f(m)|^2 |\sum_{m \leq y/p} f(m)|^2 } \Biggr) \nonumber \\
&& + O\Biggl(\frac{1}{\eta^3} \sum_{X < p \leq X^{4/3}} (\sum_{x \in \mathcal{X}} \frac{1}{\sqrt{x}} |\sum_{m \leq x/p} f(m)|)^3 \Biggr) .
\end{eqnarray}
Here the $Z_x$ are jointly normal random variables, each having mean value $\E Z_x := \tilde{\E} \Re \frac{1}{\sqrt{x}} \sum_{X < p \leq x} f(p) \sum_{m \leq x/p} f(m) = 0$ (since the $(f(p))_{p > X}$ have mean zero and are independent of the conditioning), and with variances
\begin{eqnarray}
\E Z_{x}^2 & := & \tilde{\E} \Biggl(\Re \frac{1}{\sqrt{x}} \sum_{X < p \leq x} f(p) \sum_{m \leq x/p} f(m) \Biggr)^2 \nonumber \\
& = & \tilde{\E}\Biggl(\frac{\sum_{X < p \leq x} f(p) \sum_{m \leq \frac{x}{p}} f(m) + \sum_{X < p \leq x} \overline{f(p)} \sum_{m \leq \frac{x}{p}} \overline{f(m)}}{2\sqrt{x}} \Biggr)^2 = \frac{1}{2x} \sum_{X < p \leq x} |\sum_{m \leq \frac{x}{p}} f(m)|^2 , \nonumber
\end{eqnarray}
and with covariances $\E Z_x Z_y = \frac{1}{2\sqrt{xy}} \sum_{X < p \leq X^{4/3}} \Re (\sum_{m \leq \frac{x}{p}} f(m)) (\sum_{m \leq \frac{y}{p}} \overline{f(m)})$.

We would like to show that with high probability over realisations of the $(f(p))_{p \leq X}$, the ``big Oh'' terms in \eqref{gaussapproxdisplay} will be small. To do this, we can first upper bound the sums over $\mathcal{X}$ there by sums over all points $x = X^{8/7} e^{2\pi r}$ with $0 \leq r \leq 2(\log X)/(21\pi)$. Then using H\"older's inequality and bounds for integer moments of $\sum f(m)$ (see e.g. Theorem 1.1 of Harper~\cite{harperrmfhigh}, although weaker and simpler bounds would suffice here), and the fact we sum over $\ll \log X$ values of $x$, we find the expectation of those ``big Oh'' terms is
\begin{eqnarray}
& \ll & \frac{1}{\eta^{2}} \sum_{x,y} \sqrt{\sum_{\substack{X < p \\ \leq X^{4/3}}} \frac{1}{xy} \sqrt{\E|\sum_{m \leq \frac{x}{p}} f(m)|^4 \E |\sum_{m \leq \frac{y}{p}} f(m)|^4}} + \frac{1}{\eta^3} \sum_{\substack{X < p \\ \leq X^{4/3}}} \log^{2}X \sum_{x} \frac{1}{x^{3/2}} \E|\sum_{m \leq \frac{x}{p}} f(m)|^3 \nonumber \\
& \ll & \frac{1}{\eta^{2}} \sum_{x,y} \sqrt{\sum_{\substack{X < p \\ \leq X^{4/3}}} \frac{1}{xy} \frac{x}{p} \frac{y}{p} \log X } + \frac{1}{\eta^3} \sum_{\substack{X < p \\ \leq X^{4/3}}} \log^{2}X \sum_{x} \frac{1}{x^{3/2}} (\frac{x}{p})^{3/2} \log^{1/4}X \ll \frac{\log^{9/4}X}{\eta^3 \sqrt{X}} . \nonumber
\end{eqnarray}
This is extremely small, in particular Markov's inequality now implies that with probability $\geq 1 - O(X^{-0.2})$ over all realisations of the $(f(p))_{p \leq X}$, the ``big Oh'' terms in \eqref{gaussapproxdisplay} will be $\ll \frac{\log^{9/4}X}{\eta^3 X^{0.3}} \ll \frac{1}{\eta^{3} X^{1/4}}$, say.

\vspace{12pt}
Next we want to use Normal Comparison Result 1 to understand $\p(\max_{x \in \mathcal{X}} Z_x \leq u + \eta )$ in \eqref{gaussapproxdisplay}. To do this effectively, we require more explicit estimates for the variances $\E Z_{x}^2$ and covariances $\E Z_x Z_y$. We can obtain the first of these relatively straightforwardly by combining known arguments with calculations we have already performed.

\begin{prop1}\label{propvariances}
Let $f(n)$ be a Steinhaus or Rademacher random multiplicative function. Then uniformly for all large $X$ and $X^{8/7} \leq x \leq X^{4/3}$, we have $\frac{1}{x} \sum_{X < p \leq x} |\sum_{m \leq \frac{x}{p}} f(m)|^2 \gg \frac{1}{\log X} \int_{1}^{X^{1/7}} |\sum_{m \leq z} f(m)|^2 \frac{dz}{z^2}$.

Furthermore, for any $1 \leq W \leq (\log\log X)^{1/100}$ (say) we have
$$ \p\left(\frac{1}{\log X} \int_{1}^{X^{1/7}} |\sum_{m \leq z} f(m)|^2 \frac{dz}{z^2} \geq \frac{1}{e^{2.2W} \sqrt{\log\log X}} \right) \geq 1 - O(e^{-0.1W}) . $$
\end{prop1}

\begin{proof}[Proof of Proposition \ref{propvariances}]
The first statement, which is completely deterministic, follows by an argument from the proof of Proposition 4.3 of Harper~\cite{harperrmfhigh}. Notice first that
$$ \sum_{X < p \leq x} |\sum_{m \leq x/p} f(m)|^2 \geq \frac{1}{\log x} \sum_{X < p \leq x} \log p \Biggl|\sum_{m \leq x/p} f(m) \Biggr|^2 = \frac{1}{\log x} \sum_{r \leq \frac{x}{X}} \sum_{x/(r+1) < p \leq x/r} \log p |\sum_{m \leq r} f(m)|^2 . $$
Here we have $r \leq x/X \leq x/x^{3/4} = x^{1/4}$, and so $x/r - x/(r+1) = x/(r(r+1)) \gg (x/r)^{2/3}$ and a Hoheisel-type prime number theorem in short intervals implies the above is
$$ \gg \frac{1}{\log x} \sum_{r \leq \frac{x}{X}} \Biggl(\int_{x/(r+1)}^{x/r} 1 dt \Biggr) |\sum_{m \leq r} f(m)|^2 \geq \frac{1}{\log x} \int_{X}^{x} |\sum_{m \leq x/t} f(m)|^2 dt . $$
Making the usual substitution $z=x/t$, this is all $= \frac{x}{\log x} \int_{1}^{x/X} |\sum_{m \leq z} f(m)|^2 \frac{dz}{z^2} \gg \frac{x}{\log X} \int_{1}^{X^{1/7}} |\sum_{m \leq z} f(m)|^2 \frac{dz}{z^2}$, and the claimed bound follows on dividing through by $x$.

For the second statement, we can first lower bound $\frac{1}{\log X} \int_{1}^{X^{1/7}} |\sum_{m \leq z} f(m)|^2 \frac{dz}{z^2}$ by $\frac{1}{\log X} \int_{1}^{X^{1/7}} |\sum_{\substack{m \leq z, \\ m \; \text{is} \; X \; \text{smooth}}} f(m) |^2 \frac{dz}{z^{2+84W/\log X}}$, which is
$$ \geq \frac{1}{\log X} \Biggl( \int_{1}^{\infty} \Biggl|\sum_{\substack{m \leq z, \\ X \; \text{smooth}}} f(m) \Biggr|^2 \frac{dz}{z^{2+84W/\log X}} - \frac{1}{e^{6W}} \int_{X^{1/7}}^{\infty} \Biggl|\sum_{\substack{m \leq z, \\ X \; \text{smooth}}} f(m) \Biggr|^2 \frac{dz}{z^{2+42W/\log X}} \Biggr) . $$
Using Harmonic Analysis Result 1, we can further lower bound all this by
$$ \frac{1}{2\pi \log X} \Biggl( \int_{1/3}^{1/2} \left|\frac{F(1/2 + 42W/\log X + it)}{1/2 + 42W/\log X + it}\right|^2 dt - \frac{1}{e^{6W}} \int_{-\infty}^{\infty} \left|\frac{F(1/2 + 21W/\log X + it)}{1/2 + 21W/\log X + it}\right|^2 dt \Biggr) , $$
where $F(s)$ is the Euler product corresponding to $f(m)$ on $X$-smooth numbers. Now the first integral is $\geq \int_{1/3}^{1/2} \left|F(1/2 + 42W/\log X + it) \right|^2 dt$, and Multiplicative Chaos Result 4 implies that with probability $\geq 1 - O(e^{-0.1W})$ this will be $\geq \frac{\log X}{e^{2.1W} 42 W \sqrt{\log\log X}}$. Furthermore, applying Multiplicative Chaos Result 3 with the choice $q=2/3$, say, (along with H\"older's inequality to handle the tails of the integral,) we obtain that $\E ( \int_{-\infty}^{\infty} \left|\frac{F(1/2 + 21W/\log X + it)}{1/2 + 21W/\log X + it}\right|^2 dt )^{2/3}$ is
\begin{eqnarray}
& \ll & \sum_{|N| \leq \lceil \log X \rceil} \E \Biggl( \int_{N - 1/2}^{N + 1/2} \frac{|F(\frac{1}{2} + \frac{21W}{\log X} + it)|^2}{1+N^2} \Biggr)^{2/3} + \Biggl( \E \int_{|t| > \log X} \frac{|F(\frac{1}{2} + \frac{21W}{\log X} + it)|^2}{t^2} \Biggr)^{2/3} \nonumber \\
& \ll & \sum_{|N| \leq \lceil \log X \rceil} \frac{(\log\log(|N| + 10))^{2/3}}{(1 + N^2)^{2/3}} \Biggl(\frac{\log X}{W \sqrt{\log\log X}} \Biggr)^{2/3} + 1 \ll \Biggl(\frac{\log X}{W \sqrt{\log\log X}} \Biggr)^{2/3} . \nonumber
\end{eqnarray}
Then Markov's inequality applied to the $2/3$ power implies that, with probability $\geq 1 - O(e^{-2W})$, the subtracted integral $\int_{-\infty}^{\infty} \left|\frac{F(1/2 + 21W/\log X + it)}{1/2 + 21W/\log X + it}\right|^2 dt$ will be $\leq \frac{e^{3W} \log X}{W \sqrt{\log\log X}}$.

Finally we obtain, as desired, that $\frac{1}{\log X} \int_{1}^{X^{1/7}} |\sum_{m \leq z} f(m)|^2 \frac{dz}{z^2} \geq \frac{1}{2\pi W \sqrt{\log\log X}} (\frac{1}{42 e^{2.1W}} - \frac{1}{e^{3W}}) \geq \frac{1}{e^{2.2W} \sqrt{\log\log X}}$ with probability $\geq 1 - O(e^{-0.1W}) - O(e^{-2W}) = 1 - O(e^{-0.1W})$ (assuming without loss of generality that $W$ is large enough).
\end{proof}

Over the course of the next three subsections, we shall prove the following proposition allowing us to handle the covariances.

\begin{prop2}\label{propcovariances}
Let $f(n)$ be a Steinhaus or Rademacher random multiplicative function, let $X$ be large, and for each $x = X^{8/7} e^{2\pi r}$ with $0 \leq r \leq 2(\log X)/(21\pi)$ set
$$ \mathcal{B}_{x} := \{y = X^{8/7} e^{2\pi u} : u \leq \frac{2\log X}{21\pi}, \;\;\; |\frac{1}{\sqrt{xy}} \sum_{X < p \leq X^{4/3}} (\sum_{m \leq \frac{x}{p}} f(m)) (\sum_{m \leq \frac{y}{p}} \overline{f(m)})| \geq \frac{1}{(\log\log X)^{0.6}} \} . $$
Then with probability $\geq 1 - O((\log\log X)^{-0.1})$, we have $\max_{x} \#\mathcal{B}_{x} \leq \log^{0.7}X$.
\end{prop2}

\vspace{12pt}
Assuming the truth of Proposition \ref{propcovariances}, we now finish the proof of Theorem \ref{thmlocalas} in the Steinhaus case. We assume, without loss of generality, that $W \leq (1/30)\log\log\log X$ is large. We take $u = \frac{2(\log\log X)^{1/4}}{e^{1.2W}}$ in \eqref{gaussapproxdisplay}, with $\eta$ a small {\em fixed} constant, and seek an upper bound for $\p(\max_{x \in \mathcal{X}} Z_x \leq \frac{2(\log\log X)^{1/4}}{e^{1.2W}} + \eta )$. Recall that the variances and covariances of the Gaussian $Z_x$, as well as the set $\mathcal{X}$, depend on $(f(p))_{p \leq X}$.

If we let $V := \min_{x \in \mathcal{X}} \E Z_{x}^2$ (which again depends on $(f(p))_{p \leq X}$), we certainly have
\begin{equation}\label{maxsimpeq}
\p(\max_{x \in \mathcal{X}} Z_x \leq \frac{2(\log\log X)^{1/4}}{e^{1.2W}} + \eta ) \leq \p\Biggl(\max_{x \in \mathcal{X}} \frac{Z_x}{\sqrt{\E Z_x^2}} \leq \frac{1}{\sqrt{V}} (\frac{2(\log\log X)^{1/4}}{e^{1.2W}} + \eta ) \Biggr) .
\end{equation}
By Proposition \ref{propvariances}, with probability $\geq 1 - O(e^{-0.1W})$ over all realisations of $(f(p))_{p \leq X}$ we have $V \gg \frac{1}{e^{2.2W} \sqrt{\log\log X}}$, and so $\frac{1}{\sqrt{V}} (\frac{2(\log\log X)^{1/4}}{e^{1.2W}} + \eta ) \leq \sqrt{0.5\log\log X}$, say. Furthermore, the $Z_{x}/\sqrt{\E Z_{x}^2}$ are mean zero, variance one, jointly normal random variables.

Now by Proposition \ref{propcovariances}, with probability $\geq 1 - O((\log\log X)^{-0.1})$ over all realisations of $(f(p))_{p \leq X}$ we have $\max_{x} \#\mathcal{B}_{x} \leq \log^{0.7}X$. If we have this bound, then we can greedily select a subset $\mathcal{X}' \subseteq \mathcal{X}$, of size at least $\frac{\#\mathcal{X}}{1 + \log^{0.7}X} \geq (1+o(1)) \frac{\log^{0.3}X}{11\pi}$, such that
$$ \left|\E \frac{Z_x}{\sqrt{\E Z_x^2}} \frac{Z_y}{\sqrt{\E Z_y^2}}\right| \leq \frac{(1/2) (\log\log X)^{-0.6}}{V} \ll \frac{e^{2.2W}}{(\log\log X)^{0.1}} \leq 10^{-4} \;\;\;\;\; \forall \; x,y \in \mathcal{X}' , \;\;\; x \neq y , $$
say. Here we used our lower bound $\#\mathcal{X} \geq (\log X)/(11\pi)$ (holding with probability $1 - O((\log\log X)^{-0.02})$ over realisations of $(f(p))_{p \leq X}$), along with the fact that $\E Z_x Z_y = \frac{1}{2\sqrt{xy}} \sum_{X < p \leq X^{4/3}} \Re (\sum_{m \leq \frac{x}{p}} f(m)) (\sum_{m \leq \frac{y}{p}} \overline{f(m)})$, and our upper bound assumption $W \leq (1/30)\log\log\log X$. Once again, the set $\mathcal{X}'$ depends on $(f(p))_{p \leq X}$ but is fixed having conditioned on those values. Inserting all this into Normal Comparison Result 1, noting $\sqrt{0.5\log\log X} \leq \sqrt{1.99\log(\#\mathcal{X}')}$, we see the right hand side of \eqref{maxsimpeq} will be
$$ \leq \p\left(\max_{x \in \mathcal{X}'} \frac{Z_x}{\sqrt{\E Z_x^2}} \leq \sqrt{1.99\log(\#\mathcal{X}')} \right) \ll e^{-\Theta(\log^{c}X)} + (\#\mathcal{X}')^{-1/50} \ll (\log X)^{-3/500} . $$

Collecting together the probability bounds for our various events, we have shown
\begin{eqnarray}
&& \p\Biggl(\max_{0 \leq r \leq \frac{2\log X}{21\pi}} \Re \frac{1}{\sqrt{X^{8/7} e^{2\pi r}}} \sum_{n \leq X^{8/7}e^{2\pi r}} f(n) \leq \frac{2(\log\log X)^{1/4}}{e^{1.2W}} - (\log\log X)^{0.01} \Biggr) \nonumber \\
& \leq & \p\Biggl(\max_{r \in \mathcal{R}} \Re \frac{1}{\sqrt{X^{8/7} e^{2\pi r}}} \sum_{\substack{X < p \\ \leq X^{8/7}e^{2\pi r}}} f(p) \sum_{m \leq \frac{X^{8/7}e^{2\pi r}}{p}} f(m) \leq \frac{2(\log\log X)^{1/4}}{e^{1.2W}} \Biggr) + O(\frac{1}{(\log\log X)^{0.02}}) \nonumber \\
& \ll & \frac{1}{(\log\log X)^{0.02}} + \frac{1}{X^{0.2}} + \frac{1}{e^{0.1W}} + \frac{1}{(\log\log X)^{0.1}} + \frac{1}{\eta^{3} X^{1/4}} + \frac{1}{\log^{3/500}X} \ll \frac{1}{e^{0.1W}} , \nonumber
\end{eqnarray}
which implies Theorem \ref{thmlocalas} since $\frac{2(\log\log X)^{1/4}}{e^{1.2W}} - (\log\log X)^{0.01} \geq \frac{(\log\log X)^{1/4}}{e^{1.2W}}$.
\qed

\vspace{12pt}
In the Rademacher case, the argument is identical until the calculation of the variances and covariances. One finds (since taking the real part is now redundant) that $\E Z_{x}^2 = \frac{1}{x} \sum_{X < p \leq x} (\sum_{m \leq \frac{x}{p}} f(m))^2$ and $\E Z_x Z_y = \frac{1}{\sqrt{xy}} \sum_{X < p \leq X^{4/3}} (\sum_{m \leq \frac{x}{p}} f(m)) (\sum_{m \leq \frac{y}{p}} f(m))$, i.e. differing from the Steinhaus case by a factor of 2. This is unimportant in all the subsequent calculations. The only other difference comes when taking the expectation of the ``big Oh'' terms in \eqref{gaussapproxdisplay}, where one ends up with slightly different powers of logarithms because the Rademacher moments involve different logarithmic terms than the Steinhaus ones. This is also unimportant compared with the power saving one has there, so the proof goes through in the Rademacher case as well.
\qed

\subsection{Manipulating the covariances}\label{subseccovarmanip}
Our goal now is to obtain a more tractable expression for the sums $\frac{1}{\sqrt{xy}} \sum_{X < p \leq X^{4/3}} (\sum_{m \leq \frac{x}{p}} f(m)) (\sum_{m \leq \frac{y}{p}} \overline{f(m)})$ arising in our covariances, where $X^{8/7} \leq x,y \leq X^{4/3}$. Unlike with the variances in Proposition \ref{propvariances}, we don't have positivity here so there seems no obvious direct way to relate these sums to an infinite integral average (to which we could apply some version of Parseval's identity).

Instead, using the truncated Perron formula (see e.g. Theorem 5.2 and Corollary 5.3 of Montgomery and Vaughan~\cite{mv}) we can write
$$ \frac{1}{\sqrt{x}} \sum_{m \leq \frac{x}{p}} f(m) = \frac{1}{2\pi i \sqrt{x}} \int_{1-iX^{3/4}}^{1+iX^{3/4}} F(s) \frac{(x/p)^s}{s} ds + O\left(\frac{g(x/p)}{\sqrt{x}} + \frac{\sqrt{x} X^{1/2}}{p X^{3/4}} \log X + \frac{\log X}{\sqrt{x} X^{3/4}} \right) , $$
where $F(s)$ is the finite Euler product corresponding to our random multiplicative function $f(m)$ on $X$-smooth numbers (say), and $g(x/p)$ is 1 if the distance from $x/p$ to the nearest natural number is $\leq X^{-1/2}$, and is zero otherwise. On our ranges of $p$ and $x$, note that the ``big Oh'' term is $\ll \frac{g(x/p)}{\sqrt{X}} + \frac{1}{\sqrt{p} X^{0.01}}$. Next, using Cauchy's Residue Theorem we move the line of integration and rewrite $\frac{1}{2\pi i \sqrt{x}} \int_{1-iX^{3/4}}^{1+iX^{3/4}} F(s) \frac{(x/p)^s}{s} ds$ as
\begin{eqnarray}
&& \frac{1}{2\pi i \sqrt{x}} \int_{1/2-iX^{3/4}}^{1/2+iX^{3/4}} F(s) \frac{x^s}{s p^s} + O(\frac{1}{\sqrt{x} X^{3/4}} \int_{1/2}^{1} (|F(\sigma + iX^{3/4})| + |F(\sigma - iX^{3/4})|) (\frac{x}{p})^{\sigma} d\sigma ) \nonumber \\
& = & \frac{1}{2\pi i \sqrt{x}} \int_{1/2-iX^{3/4}}^{1/2+iX^{3/4}} F(s) \frac{x^s}{s p^s} + O(\frac{1}{X^{3/4}} (\frac{1}{\sqrt{x}} + \frac{\sqrt{x}}{p}) \int_{1/2}^{1} (|F(\sigma + iX^{3/4})| + |F(\sigma - iX^{3/4})|) ) , \nonumber
\end{eqnarray}
where the second line follows after upper bounding $(x/p)^{\sigma}$ by $1 + (x/p)$. The integral remaining inside the ``big Oh'' term, which doesn't involve $x$ or $p$, has expectation $\ll \int_{1/2}^{1} \sqrt{\E|F(\sigma + iX^{3/4})|^2} d\sigma \ll \sqrt{\log X}$, so by Markov's inequality it will be $\leq \log X$ with probability $\geq 1 - O(\log^{-1/2}X)$. If this event occurs, then this entire ``big Oh'' term will be $\ll \frac{\log X}{X^{3/4}} (\frac{1}{\sqrt{x}} + \frac{\sqrt{x}}{p}) \ll \frac{1}{\sqrt{p} X^{0.01}}$ for {\em all} $x$ and $p$.

Inserting the above expression for the sums up to $x/p$ and $y/p$ and summing over $p$, we deduce that (with probability $\geq 1 - O(\log^{-1/2}X)$) we have
\begin{eqnarray}\label{afterperrondisplay}
&& \frac{1}{\sqrt{xy}} \sum_{X < p \leq X^{4/3}} (\sum_{m \leq \frac{x}{p}} f(m)) (\sum_{m \leq \frac{y}{p}} \overline{f(m)}) \\
& = & \frac{1}{4\pi^2} \int_{-X^{3/4}}^{X^{3/4}} \int_{-X^{3/4}}^{X^{3/4}} \frac{F(1/2+iv) x^{iv}}{1/2 + iv} \frac{\overline{F(1/2+it)} y^{-it}}{1/2 - it} \sum_{X < p \leq X^{4/3}} \frac{1}{p^{1+i(v - t)}} dv dt + \nonumber \\
&& + O\Biggl(\sum_{X < p \leq X^{4/3}} \Biggl((\frac{g(\frac{x}{p})+g(\frac{y}{p})}{\sqrt{Xp}} + \frac{1}{p X^{0.01}}) \int_{-X^{3/4}}^{X^{3/4}} \frac{|F(\frac{1}{2}+it)|}{|1/2 + it|} dt + \frac{g(\frac{x}{p})+g(\frac{y}{p})}{X} + \frac{1}{p X^{0.02}} \Biggr) \Biggr) \nonumber
\end{eqnarray}
for all $X^{8/7} \leq x,y \leq X^{4/3}$. Again, the integral inside the ``big Oh'' term doesn't depend on $x$ or $y$, and has expectation $\leq \int_{-X^{3/4}}^{X^{3/4}} \frac{\sqrt{\E|F(1/2+it)|^2}}{|1/2 + it|} dt \ll \log^{3/2}X$, so by Markov's inequality it will be $\leq \log^{2}X$ with probability $\geq 1 - O(\log^{-1/2}X)$. If this event occurs, then the ``big Oh'' term in \eqref{afterperrondisplay} will be $\ll \frac{\log^{2}X}{X} \sum_{X < p \leq X^{4/3}} (g(x/p) + g(y/p)) + \frac{\log^{2}X}{X^{0.01}} \ll \frac{\log^{2}X}{X} \frac{X^{4/3}}{X^{1/2}} + \frac{\log^{2}X}{X^{0.01}} \ll \frac{\log^{2}X}{X^{0.01}}$ (which is negligibly small) for all $X^{8/7} \leq x,y \leq X^{4/3}$.

\vspace{12pt}
Next we want to clean up \eqref{afterperrondisplay} further by restricting the ranges of integration. By symmetry between $v, t$, the part of the double integral where $\max\{|v|,|t|\} \geq \log^{2}X$ is
$$ \ll \int_{-X^{3/4}}^{X^{3/4}} \frac{|F(1/2+it)|}{1+|t|} \int_{\max\{|t|, \log^{2}X\} \leq |v| \leq X^{3/4}} \frac{|F(1/2+iv)|}{1+|v|} |\sum_{X < p \leq X^{4/3}} \frac{1}{p^{1+i(v - t)}}| dv dt . $$
As in all our earlier calculations, the expectation of this expression is crudely $\ll \log X \int_{-X^{3/4}}^{X^{3/4}} \frac{1}{1+|t|} \int_{\max\{|t|, \log^{2}X\} \leq |v| \leq X^{3/4}} \frac{1}{1+|v|} |\sum_{X < p \leq X^{4/3}} \frac{1}{p^{1+i(v - t)}}| dv dt$. Furthermore, applying the Cauchy--Schwarz inequality followed by Number Theory Result 1 (with $a_n$ chosen to be $n^{it}/\log n$ when $n$ is a prime between $X$ and $X^{4/3}$, and chosen to be zero otherwise) and standard Chebychev-type prime number estimates, the inner integral over $v$ here is at most
$$ \sqrt{\int_{|v| \geq \max\{|t|, \log^{2}X\}} \frac{dv}{(1+|v|)^2}} \sqrt{\int_{-X^{3/4}}^{X^{3/4}} |\sum_{X < p \leq X^{4/3}} \frac{1}{p^{1+i(v - t)}}|^2 dv} \ll \sqrt{\frac{1}{\max\{|t|, \log^{2}X\}} \frac{1}{\log X}} . $$
Integrating over $t$, we obtain $\log X \int_{-\log^{2}X}^{\log^{2}X} \frac{dt}{(1+|t|)\log^{3/2}X} + \log X \int_{|t| \geq \log^{2}X} \frac{dt}{|t|^{3/2} \log^{1/2}X} \ll \frac{\log\log X}{\sqrt{\log X}}$. So by Markov's inequality, with probability $\geq 1 - O(\log^{-1/4}X)$ we can restrict the range of integration in \eqref{afterperrondisplay} to $\max\{|v|,|t|\} \leq \log^{2}X$, at the cost of a negligible error term $O(\frac{\log\log X}{\log^{1/4}X})$.

Finally we will show that with high probability, we can restrict the range of integration in \eqref{afterperrondisplay} to the ``near diagonal'' portion where $|v - t| \leq \frac{(\log\log X)^{100}}{\log X}$, say. Using Number Theory Result 2, the complementary portion of the double integral contributes
$$ \ll \int \int_{\substack{|v|, |t| \leq \log^{2}X, \\ |v - t| > \frac{(\log\log X)^{100}}{\log X}}}  \frac{|F(1/2+iv)|}{1+ |v|} \frac{|F(1/2+it)|}{1 + |t|} \frac{1}{|v - t| \log X} dv dt . $$
In the Steinhaus case, the second part of Euler Product Result 1 implies this has expectation
$$ \ll \frac{1}{\sqrt{\log X}} \int \int_{\substack{|v|, |t| \leq \log^{2}X, \\ |v - t| > \frac{(\log\log X)^{100}}{\log X}}}  \frac{1}{1+ |v|} \frac{1}{1 + |t|} \frac{1}{|v - t|} \left(\frac{\textbf{1}_{|v - t| \leq 1}}{\sqrt{|v - t|}} + \log(2+|v - t|) \right) dv dt . $$
If $|v - t| \leq 1$ then $1+|v| \asymp 1+|t|$, so the corresponding contribution to the integral is $\ll \frac{1}{\sqrt{\log X}} \int_{-\log^{2}X}^{\log^{2}X} \frac{1}{(1+ |v|)^2} \int_{\frac{(\log\log X)^{100}}{\log X}}^{1} \frac{dh}{h^{3/2}} dv \ll \frac{1}{(\log\log X)^{50}}$. Crudely, the double integral of $\log(2+|v - t|)$ is $\ll \frac{\log\log X}{\sqrt{\log X}} \int_{-\log^{2}X}^{\log^{2}X} \frac{1}{1+ |v|} \int_{\frac{(\log\log X)^{100}}{\log X}}^{2\log^{2}X} \frac{dh}{h} dv \ll \frac{(\log\log X)^3}{\sqrt{\log X}}$. In the Rademacher case, the expectation calculations are very similar except that the second part of Euler Product Result 2 produces an additional term $\frac{\textbf{1}_{|v + t| \leq 1}}{\sqrt{|v + t|}}$ inside the bracket in the double integral. Again, if $|v + t| \leq 1$ then we have $1+|v| \asymp 1+|t|$, so the extra contribution to the integral is
$$ \ll \frac{1}{\sqrt{\log X}} \int_{-\log^{2}X}^{\log^{2}X} \frac{1}{(1+ |v|)^2} \int_{0}^{1} \min\{\frac{1}{|2|v| - h|}, \frac{\log X}{(\log\log X)^{100}}\} \frac{dh}{\sqrt{h}} dv . $$
Switching the order of the integrals, this is readily seen to be $\ll \frac{\log\log X}{\sqrt{\log X}}$. So finally in both the Steinhaus and Rademacher cases, Markov's inequality yields that with probability $\geq 1 - O((\log\log X)^{-25})$, we can restrict to the near diagonal portion at the cost of a negligible error term $O((\log\log X)^{-25})$.

\subsection{Inserting a strong barrier condition}\label{subsecstrongbarrier}
Recall that in view of the calculations performed in section \ref{subseccovarmanip}, with probability $\geq 1 - O((\log\log X)^{-25})$ the covariance expression $\frac{1}{\sqrt{xy}} \sum_{X < p \leq X^{4/3}} (\sum_{m \leq \frac{x}{p}} f(m)) (\sum_{m \leq \frac{y}{p}} \overline{f(m)})$ is
$$ = \frac{1}{4\pi^2} \int \int_{\substack{|v|, |t| \leq \log^{2}X , \\ |v - t| \leq \frac{(\log\log X)^{100}}{\log X}}} \frac{F(\frac{1}{2}+iv) x^{iv}}{1/2 + iv} \frac{\overline{F(\frac{1}{2}+it)} y^{-it}}{1/2 - it} \sum_{X < p \leq X^{4/3}} \frac{1}{p^{1+i(v - t)}} + O((\log\log X)^{-25}) $$
for all $x,y$ as in Proposition \ref{propcovariances}. For orientation, we note that since $v, t$ are so close and since $\sum_{X < p \leq X^{4/3}} \frac{1}{p^{1+i(v - t)}}$ also decays in size quite rapidly once $|v - t| > 1/\log X$, this expression is (heuristically at least) $\approx \frac{1}{\log X} \int_{-\log^{2}X}^{\log^{2}X} \frac{|F(1/2+it)|^2}{|1/2+it|^2} (x/y)^{it} dt$, similarly as in the discussion in the Introduction.

Our next objective is to edit the double integral even more, in preparation for our final proof of Proposition \ref{propcovariances}. Unlike in section \ref{subseccovarmanip}, which was mostly of a general complex analytic flavour and where we had much room to spare with most estimates, these final edits will invoke multiplicative chaos machinery from section \ref{subsecmultchaos} and be more delicate.

Applying Multiplicative Chaos Result 1 with the choice $e^{W} = (\log\log X)^3$, say, we find that for any given $N \in \R$ we have
$$ |F_{j}(1/2 + it(j))| \leq \frac{\log X}{e^j} (\log\log X - j)^2 (\log\log X)^3 \leq \frac{\log X}{e^j} (\log\log X)^5 $$
for all $0 \leq j \leq \log\log X - 1$ and all $|t-N| \leq 1/2$, with probability $\geq 1 - O((\log\log X)^{-6})$. Recall here from section \ref{sectools} that $F_{j}(s)$ denotes the partial random Euler product over $X^{e^{-j}}$-smooth numbers, and $t(j)$ denotes the sequence of increasingly coarse approximations to $t$. In particular, if we let $\mathcal{D}(t)$ denote the event that $|F_{j}(1/2 + it(j))| \leq \frac{\log X}{e^j} (\log\log X)^5$ for all $0 \leq j \leq \log\log X - 1$, and let $\mathcal{D}$ be the event that $\mathcal{D}(t)$ occurs for all $|t| \leq (\log\log X)^2$ (say), then the union bound implies that $\mathcal{D}$ holds with probability $\geq 1 - O((\log\log X)^{-4})$. Furthermore, the contribution to our covariance double integral from those $v, t$ with $\max\{|v|,|t|\} \geq (\log\log X)^2$ is always
$$ \ll \int_{(\log\log x)^2 \leq |v| \leq \log^{2}X} \int_{|v - t| \leq \frac{(\log\log X)^{100}}{\log X}} \frac{|F(\frac{1}{2}+iv)|}{|1/2 + iv|^2} |F(\frac{1}{2}+it)| |\sum_{X < p \leq X^{4/3}} \frac{1}{p^{1+i(v - t)}}| dt dv , $$
which crudely (and using Number Theory Result 2 to bound the prime sum) has expectation $\ll \log X \cdot \frac{1}{(\log\log X)^2} \cdot \int_{|h| \leq \frac{(\log\log X)^{100}}{\log X}} |\sum_{X < p \leq X^{4/3}} \frac{1}{p^{1+ih}}| dh \ll \frac{\log\log\log X}{(\log\log X)^2}$. Applying Markov's inequality to this and using our bound for $\p(\mathcal{D} \; \text{holds})$, we find that with probability $\geq 1 - O((\log\log X)^{-1})$, the sum $\frac{1}{\sqrt{xy}} \sum_{X < p \leq X^{4/3}} (\sum_{m \leq \frac{x}{p}} f(m)) (\sum_{m \leq \frac{y}{p}} \overline{f(m)})$ is
$$ \frac{1}{4\pi^2} \int \int_{\substack{|v|, |t| \leq (\log\log X)^{2} , \\ |v - t| \leq \frac{(\log\log X)^{100}}{\log X}}} \frac{\textbf{1}_{\mathcal{D}(v)} F(\frac{1}{2}+iv) x^{iv}}{1/2 + iv} \frac{\textbf{1}_{\mathcal{D}(t)} \overline{F(\frac{1}{2}+it)} y^{-it}}{1/2 - it} \sum_{\substack{X < p \\ \leq X^{4/3}}} \frac{1}{p^{1+i(v - t)}} + O(\frac{\log\log\log X}{\log\log X}) $$
for all $x,y$ as in Proposition \ref{propcovariances}. Here $\textbf{1}$ as usual denotes the indicator function.

Finally we use Multiplicative Chaos Result 2 to strengthen the barrier condition $\mathcal{D}(t)$. To ensure this result is applicable (and to guarantee that $v, t$ always have the same sign, which will slightly simplify matters in the next section), we note that the portion of the double integral where $\min\{|v|,|t|\} \leq \frac{1}{(\log\log X)^2}$ (say) has expectation $\ll \log X \cdot \frac{1}{(\log\log X)^2} \cdot \int_{|h| \leq \frac{(\log\log X)^{100}}{\log X}} |\sum_{X < p \leq X^{4/3}} \frac{1}{p^{1+ih}}| dh \ll \frac{\log\log\log X}{(\log\log X)^2}$, so will be $\ll \frac{\log\log\log X}{\log\log X}$ with probability $\geq 1 - O((\log\log X)^{-1})$ by Markov's inequality. Then if $\mathcal{A}(t)$ denotes the event that
$$ |F_{j}(\frac{1}{2} + it)| \leq \frac{\log X}{e^j (\log\log X)^{1000}} \;\; \forall \; j \leq 0.99\log\log X , \;\; \text{and} \;\; |F_{j}(\frac{1}{2} + it)| \leq \frac{\log X}{e^j} (\log\log X)^6 \;\; \forall \; j , $$
with probability $\geq 1 - O((\log\log X)^{-1})$ we have
\begin{eqnarray}\label{afterbarrierdisplay}
&& \frac{1}{\sqrt{xy}} \sum_{X < p \leq X^{4/3}} (\sum_{m \leq \frac{x}{p}} f(m)) (\sum_{m \leq \frac{y}{p}} \overline{f(m)}) \\
& = & \frac{1}{4\pi^2} \int \int_{\substack{|v|, |t| \leq (\log\log X)^{2} , \\ |v|, |t| \geq (\log\log X)^{-2} , \\ |v - t| \leq \frac{(\log\log X)^{100}}{\log X}}} \frac{\textbf{1}_{\mathcal{A}(v)} F(\frac{1}{2}+iv) x^{iv}}{1/2 + iv} \frac{\textbf{1}_{\mathcal{A}(t)} \overline{F(\frac{1}{2}+it)} y^{-it}}{1/2 - it} \sum_{\substack{X < p \\ \leq X^{4/3}}} \frac{1}{p^{1+i(v - t)}} + \nonumber \\
&& + O\Biggl( \int \int \frac{\textbf{1}_{\mathcal{D}(v)} \textbf{1}_{\mathcal{A}(v) \; \text{fails}} |F(\frac{1}{2}+iv)|}{|1/2 + iv|} \frac{|F(\frac{1}{2}+it)|}{|1/2 - it|} |\sum_{\substack{X < p \\ \leq X^{4/3}}} \frac{1}{p^{1+i(v - t)}}| + \frac{\log\log\log X}{\log\log X} \Biggr) \nonumber
\end{eqnarray}
for all $x,y$ as in Proposition \ref{propcovariances}. Here the range of integration inside the ``big Oh'' term is as in the main term. Again using the Cauchy--Schwarz inequality, and the fact that $\int_{|v - t| \leq \frac{(\log\log X)^{100}}{\log X}} |\sum_{\substack{X < p \\ \leq X^{4/3}}} \frac{1}{p^{1+i(v - t)}}| dt \ll \frac{\log\log\log X}{\log X}$ and $\int_{|v - t| \leq \frac{(\log\log X)^{100}}{\log X}} |\sum_{\substack{X < p \\ \leq X^{4/3}}} \frac{1}{p^{1+i(v - t)}}| dv \ll \frac{\log\log\log X}{\log X}$ (which follows by bounding the prime number sums with Number Theory Result 2), we deduce that ``big Oh'' integral is
$$ \ll \frac{\log\log\log X}{\log X} \sqrt{\int_{(\log\log X)^{-2} \leq |v| \leq (\log\log X)^{2}} \frac{\textbf{1}_{\mathcal{D}(v)} \textbf{1}_{\mathcal{A}(v) \; \text{fails}} |F(\frac{1}{2}+iv)|^2}{|1/2 + iv|^2}} \sqrt{\int_{- (\log\log X)^{2}}^{(\log\log X)^2} \frac{|F(\frac{1}{2}+it)|^2}{|1/2 - it|^2}} . $$
By Multiplicative Chaos Result 2 and linearity of expectation we have
$$ \E \int_{(\log\log X)^{-2} \leq |v| \leq (\log\log X)^{2}} \frac{\textbf{1}_{\mathcal{D}(v)} \textbf{1}_{\mathcal{A}(v) \; \text{fails}} |F(\frac{1}{2}+iv)|^2}{|1/2 + iv|^2} \ll \log X \frac{(\log\log\log X)^7}{\log\log X} , $$
so Markov's inequality implies that $\int_{(\log\log X)^{-2} \leq |v| \leq (\log\log X)^{2}} \frac{\textbf{1}_{\mathcal{D}(v)} \textbf{1}_{\mathcal{A}(v) \; \text{fails}} |F(\frac{1}{2}+iv)|^2}{|1/2 + iv|^2}$ is $\leq \log X \frac{(\log\log\log X)^7}{(\log\log X)^{0.9}}$ with probability $\geq 1 - O((\log\log X)^{-0.1})$. Also, Multiplicative Chaos Result 3 implies $\E\left(\int_{- (\log\log X)^{2}}^{(\log\log X)^2} \frac{|F(\frac{1}{2}+it)|^2}{|1/2 - it|^2}\right)^{0.9} \leq \sum_{|N| \leq (\log\log X)^2 + 1} \E\left(\int_{N-1/2}^{N+1/2} \frac{|F(\frac{1}{2}+it)|^2}{|1/2 - it|^2}\right)^{0.9} \ll \left(\frac{\log X}{\sqrt{\log\log X}}\right)^{0.9}$, and so $\int_{- (\log\log X)^{2}}^{(\log\log X)^2} \frac{|F(\frac{1}{2}+it)|^2}{|1/2 - it|^2} \leq \frac{\log X}{(\log\log X)^{0.35}}$ with probability $\geq 1 - O((\log\log X)^{-0.1})$.

Overall we find that with probability $\geq 1 - O((\log\log X)^{-0.1})$, the ``big Oh'' term in \eqref{afterbarrierdisplay} is $\ll \frac{\log\log\log X}{\log X} \sqrt{\log X \frac{(\log\log\log X)^7}{(\log\log X)^{0.9}}} \sqrt{\frac{\log X}{(\log\log X)^{0.35}}} + \frac{\log\log\log X}{\log\log X} \ll \frac{1}{(\log\log X)^{0.61}}$, say.

\subsection{Proof of Proposition \ref{propcovariances}}\label{subsecfinalexp}
Using display \eqref{afterbarrierdisplay} and the subsequent discussion, the proof of Proposition \ref{propcovariances} is reduced to proving an analogous statement for the double integrals $\int \int_{\Delta} \frac{\textbf{1}_{\mathcal{A}(v)} F(\frac{1}{2}+iv) x^{iv}}{1/2 + iv} \frac{\textbf{1}_{\mathcal{A}(t)} \overline{F(\frac{1}{2}+it)} y^{-it}}{1/2 - it} \sum_{X < p \leq X^{4/3}} \frac{1}{p^{1+i(v - t)}} dv dt$, where
$$ \Delta := \{(v,t) : (\log\log X)^{-2} \leq |v|, |t| \leq (\log\log X)^{2} , \;\;\; |v - t| \leq \frac{(\log\log X)^{100}}{\log X} \} . $$ 
In fact, for these integrals we will be able to prove something much stronger. However, due to the error terms already accumulated in proving \eqref{afterbarrierdisplay}, this does not lead to a stronger statement of Proposition \ref{propcovariances}.

\begin{prop3}\label{proplargecovar}
For any large $X$ and any natural number $k \leq \log^{1/4}X$ (say), and with the event $\mathcal{A}(t)$ and the set $\Delta$ defined as above, the following is true. With probability $\geq 1 - (\log\log X)^{-2k}$ we have
\begin{eqnarray}
&& \max_{x = X^{8/7} e^{2\pi r}} \sum_{y = X^{8/7} e^{2\pi u}} \left| \int \int_{\Delta} \frac{\textbf{1}_{\mathcal{A}(v)} F(\frac{1}{2}+iv) x^{iv}}{1/2 + iv} \frac{\textbf{1}_{\mathcal{A}(t)} \overline{F(\frac{1}{2}+it)} y^{-it}}{1/2 - it} \sum_{X < p \leq X^{4/3}} \frac{1}{p^{1+i(v - t)}} dv dt \right|^{2k} \nonumber \\
& \ll & \log X \Biggl(\frac{1}{\log^{1/3}X} \sum_{|N| \leq (\log\log X)^2 + 1} \frac{1}{N^2 + 1} \Biggl( \frac{C \log\log\log X}{\log X} \int_{N \leq t \leq N+1} |F(\frac{1}{2} + it)|^{2} dt \Biggr)^{2k} + \nonumber \\
&& + \frac{1}{\log^{1/3}X} (k^{29} (\log\log X)^{17})^{2k} + (\frac{k^{29}}{(\log\log X)^{989}})^{2k} \Biggr) , \nonumber
\end{eqnarray}
where the maximum is over all $0 \leq r \leq (2\log X)/(21\pi)$ and the sum is over $0 \leq u \leq (2\log X)/(21\pi)$, and where $C$ is a certain absolute constant. 
\end{prop3}

\begin{proof}[Proof of Proposition \ref{propcovariances}, assuming Proposition \ref{proplargecovar}]
Using linearity of expectation and Euler Product Results 1 and 2, for all $N$ we have $\E \int_{N \leq t \leq N+1} |F(1/2 + it)|^{2} dt \ll \log X$. So by Markov's inequality, for any given $N$ we have $\int_{N \leq t \leq N+1} |F(1/2 + it)|^{2} dt \leq \log X (\log\log X)^3$ (say) with probability $\geq 1 - O((\log\log X)^{-3})$, and by the union bound this holds for all $|N| \leq (\log\log X)^2 + 1$ simultaneously with probability $\geq 1 - O((\log\log X)^{-1})$. Thus with probability $\geq 1 - O((\log\log X)^{-1})$, the left hand side in Proposition \ref{proplargecovar} will be
$$ \ll \log X \Biggl(\frac{1}{\log^{1/3}X} (k^{29} (\log\log X)^{17})^{2k} + (\frac{k^{29}}{(\log\log X)^{989}})^{2k} \Biggr) , $$
and so the maximal number of $y$ for which the double integral is $\geq \frac{1}{\log\log X}$ will be
$$ \ll \log X \Biggl(\frac{1}{\log^{1/3}X} (k^{29} (\log\log X)^{18})^{2k} + (\frac{k^{29}}{(\log\log X)^{988}})^{2k} \Biggr) . $$

Choosing $k = \lfloor \frac{\log\log X}{4000\log\log\log X} \rfloor$ to roughly balance the size of the terms in this bound, we conclude that with probability $\geq 1 - O((\log\log X)^{-1})$ the maximal number of $y$ with the double integral $\geq \frac{1}{\log\log X}$ will be $\leq \log^{0.7}X$, say. Combining this with \eqref{afterbarrierdisplay} and the subsequent discussion, Proposition \ref{propcovariances} follows.
\end{proof}

Note that to obtain a strong bound for $\#\mathcal{B}_{x}$ in Proposition \ref{propcovariances} (saving a power of $\log X$, as required to obtain large sets $\mathcal{X}'$ in the subsequent deduction of Theorem \ref{thmlocalas}), it was crucial to take a high power $2k$ in Proposition \ref{proplargecovar}. As the reader will see in the proof of Proposition \ref{proplargecovar}, this high power has the effect of boosting the saving produced by the barrier conditions $\mathcal{A}(v), \mathcal{A}(t)$ relative to the number of $y$ values being summed over.

\begin{proof}[Proof of Proposition \ref{proplargecovar}]
Firstly, if we let $\Delta_N$ denote the subset of $\Delta$ for which $N \leq \max\{t,v\} \leq N+1$, then using H\"older's inequality we can upper bound the $2k$-th power in Proposition \ref{proplargecovar} by
\begin{eqnarray}
&& \left( \sum_{N} \frac{1}{N^2 + 1} \left| \int \int_{\Delta_N} (N^2 + 1) \frac{\textbf{1}_{\mathcal{A}(v)} F(\frac{1}{2}+iv) x^{iv}}{1/2 + iv} \frac{\textbf{1}_{\mathcal{A}(t)} \overline{F(\frac{1}{2}+it)} y^{-it}}{1/2 - it} \sum_{X < p \leq X^{4/3}} \frac{1}{p^{1+i(v - t)}} \right| \right)^{2k} \nonumber \\
& \leq & C^{2k} \sum_{N} \frac{1}{N^2 + 1} \left| \int \int_{\Delta_N} (N^2 + 1) \frac{\textbf{1}_{\mathcal{A}(v)} F(\frac{1}{2}+iv) x^{iv}}{1/2 + iv} \frac{\textbf{1}_{\mathcal{A}(t)} \overline{F(\frac{1}{2}+it)} y^{-it}}{1/2 - it} \sum_{X < p \leq X^{4/3}} \frac{1}{p^{1+i(v - t)}} \right|^{2k} , \nonumber
\end{eqnarray}
where the sums are over all $-(\log\log X)^2 - 1 \leq N \leq (\log\log X)^2$, and $C$ is an absolute constant (actually one can take $C = \sum_{N=-\infty}^{\infty} \frac{1}{N^2 + 1}$ at this point). We will mostly concentrate on bounding the contribution from $N=0$. The treatment of all other $N$ is very similar, noting that the factor $N^2 + 1$ inside the absolute values is the same size, up to a constant factor, as the denominator $(1/2 + iv)(1/2 - it)$. 

Observe that for any $t_1, ..., t_{2k} \in \R$ we have
$$ \sum_{0 \leq u \leq (2\log X)/(21\pi)} e^{2\pi iu(t_1 + ... + t_k - t_{k+1} - ... - t_{2k})} \ll \min\{\log X , \frac{1}{|| t_1 + ... + t_k - t_{k+1} - ... - t_{2k} ||}\} , $$
where $||\cdot ||$ denotes distance to the nearest integer. So expanding the $2k$-th power (in the $N=0$ case) and then performing the outer summation over $y$, we obtain a bound
\begin{eqnarray}
& \ll & 4^{2k} \int_{0}^{1} ... \int_{0}^{1} \int_{|v_1-t_1| \leq \frac{(\log\log X)^{100}}{\log X}} ... \int_{|v_{2k}-t_{2k}| \leq \frac{(\log\log X)^{100}}{\log X}} \prod_{j=1}^{2k} \textbf{1}_{\mathcal{A}(t_j)} |F(\frac{1}{2} + it_j)| \textbf{1}_{\mathcal{A}(v_j)} |F(\frac{1}{2} + iv_j)| \nonumber \\
&& \cdot \prod_{j=1}^{2k} |\sum_{X < p \leq X^{4/3}} \frac{1}{p^{1+i(v_j - t_j)}}| \cdot \min\{\log X , \frac{1}{|| t_1 + ... + t_k - t_{k+1} - ... - t_{2k} ||}\} dv_{1} ... dv_{2k} dt_{1} ... dt_{2k} , \nonumber
\end{eqnarray}
where the factor $4^{2k}$ bounds the contribution from all the denominators $\frac{1}{|1/2 + iv|} \frac{1}{|1/2 - it|}$. Note that the dependence on $x$ has already disappeared here, so the maximum over $x$ in Proposition \ref{proplargecovar} becomes irrelevant.

Now we take a fairly crude approach to simplify this complicated expression. In the part of the multiple integral where $|| t_1 + ... + t_k - t_{k+1} - ... - t_{2k} || \geq 1/\log^{2/3}X$, say, we shall upper bound the minimum by $\log^{2/3}X$. Also discarding the indicator functions $\textbf{1}_{\mathcal{A}(t_j)}, \textbf{1}_{\mathcal{A}(v_j)}$, and using the fact that $\prod_{j=1}^{2k} |F(1/2 + it_j)| |F(1/2 + iv_j)| \leq \prod_{j=1}^{2k} |F(1/2 + it_j)|^2 + \prod_{j=1}^{2k} |F(1/2 + iv_j)|^2$, and using the symmetry between the $t$ variables and $v$ variables, we can bound this piece of the integral by
$$ \ll 4^{2k} \log^{2/3}X \Biggl( \int_{0}^{1} \int_{|v-t| \leq \frac{(\log\log X)^{100}}{\log X}} |F(1/2 + it)|^{2} |\sum_{X < p \leq X^{4/3}} \frac{1}{p^{1+i(v - t)}}| dv dt \Biggr)^{2k} . $$
By Number Theory Result 2 we have $|\sum_{X < p \leq X^{4/3}} \frac{1}{p^{1+i(v - t)}}| \ll \min\{1,\frac{1}{|v-t|\log X}\}$, so inserting this and integrating we can bound the above piece of the multiple integral by $\ll \log^{2/3}X \Biggl( \frac{C \log\log\log X}{\log X} \int_{0}^{1} |F(1/2 + it)|^{2} dt \Biggr)^{2k}$. Summing the contributions of this kind for all $N$ produces the first term in the Proposition \ref{proplargecovar} bound.

It remains to work with the piece of the integral where $|| t_1 + ... + t_k - t_{k+1} - ... - t_{2k} || < 1/\log^{2/3}X$. On this portion we shall upper bound $\min\{\log X , \frac{1}{|| t_1 + ... + t_k - t_{k+1} - ... - t_{2k} ||}\}$ trivially by $\log X$, but shall retain the information that we are only integrating over a small subset of the range of integration. Since $-k \leq t_1 + ... + t_k - t_{k+1} - ... - t_{2k} \leq k$, on this part of the integral we must have $|(t_1 + ... + t_k - t_{k+1} - ... - t_{2k}) - m| < 1/\log^{2/3}X$ for some integer $-k \leq m \leq k$. Furthermore, since the $v_i$ must be very close to the $t_i$ we also have $|(v_1 + ... + v_k - v_{k+1} - ... - v_{2k}) - m| < \frac{1}{\log^{2/3}X} + \frac{2k(\log\log X)^{100}}{\log X} \leq \frac{2}{\log^{2/3}X}$. So upper bounding $\prod_{j=1}^{2k} \textbf{1}_{\mathcal{A}(t_j)} |F(1/2 + it_j)| \textbf{1}_{\mathcal{A}(v_j)} |F(1/2 + iv_j)|$ by $\prod_{j=1}^{2k} \textbf{1}_{\mathcal{A}(t_j)} |F(1/2 + it_j)|^2 + \prod_{j=1}^{2k} \textbf{1}_{\mathcal{A}(v_j)} |F(1/2 + iv_j)|^2$ similarly as before, and applying the bound $|\sum_{X < p \leq X^{4/3}} \frac{1}{p^{1+i(v - t)}}| \ll \min\{1,\frac{1}{|v-t|\log X}\}$, and using the symmetry between $v$ and $t$, we can upper bound this portion of the integral by
\begin{equation}\label{integralsubeq}
\log X \left(\frac{C\log\log\log X}{\log X}\right)^{2k} \sum_{|m| \leq k} \int_{I_m} \prod_{j=1}^{2k} \textbf{1}_{\mathcal{A}(t_j)} |F(1/2+it_j)|^2 dt_1 ... dt_{2k} ,
\end{equation}
where (for $N=0$) $I_m$ is the subset of $[0,1]^{2k}$ where $|(t_1 + ... + t_k - t_{k+1} - ... - t_{2k}) - m| \leq \frac{2}{\log^{2/3}X}$.

\vspace{12pt}
Since the behaviour of $\prod_{j=1}^{2k} \textbf{1}_{\mathcal{A}(t_j)} |F(1/2+it_j)|^2$ depends on the distance between the points $t_j$ (see e.g. Euler Product Results 1 and 2), we shall rewrite each integral $\int_{I_m}$ in \eqref{integralsubeq} to make that information easier to extract. Thus we have
$$ \int_{I_m} \prod_{j=1}^{2k} \textbf{1}_{\mathcal{A}(t_j)} |F(1/2+it_j)|^2 dt_1 ... dt_{2k} = \sum_{\sigma \in \text{Sym}(2k)} \int_{I_{\sigma,m}} \prod_{j=1}^{2k} \textbf{1}_{\mathcal{A}(t_j)} |F(1/2+it_j)|^2 dt_1 ... dt_{2k} , $$
where $I_{\sigma,m}$ denotes the subset of $I_m$ for which $t_{\sigma(1)} < t_{\sigma(2)} < ... < t_{\sigma(2k)}$. And for each $\sigma \in \text{Sym}(2k)$, we can further rewrite
$$ \int_{I_{\sigma,m}} \prod_{j=1}^{2k} \textbf{1}_{\mathcal{A}(t_j)} |F(\frac{1}{2}+it_j)|^2 dt_1 ... dt_{2k} = \sum_{h_1 = 0}^{\lfloor \log\log X \rfloor} ... \sum_{h_{2k-1} = 0}^{\lfloor \log\log X \rfloor} \int_{I_{\sigma,\textbf{h},m}} \prod_{j=1}^{2k} \textbf{1}_{\mathcal{A}(t_j)} |F(\frac{1}{2}+it_j)|^2 dt_1 ... dt_{2k} , $$
where $I_{\sigma,\textbf{h},m}$ denotes the subset of $I_{\sigma,m}$ for which $\frac{e^{h_i}}{\log X} \leq t_{\sigma(i+1)} - t_{\sigma(i)} \leq \frac{e^{h_i + 1}}{\log X}$ for all $1 \leq i \leq 2k-1$ (with the convention that, whenever $h_i = 0$, the condition is instead simply that $t_{\sigma(i+1)} - t_{\sigma(i)} \leq \frac{e}{\log X}$).

Next, for ease of writing let $\alpha(h)$ equal $\frac{1}{(\log\log X)^{1000}}$ if $0 \leq h \leq 0.99\log\log X$, and equal $(\log\log X)^6$ if $h > 0.99\log\log X$. Then applying the definition of $\mathcal{A}(t)$ from section \ref{subsecstrongbarrier}, we can upper bound each integral $\int_{I_{\sigma,m}}$ by
$$ \sum_{h_1 = 0}^{\lfloor \log\log X \rfloor} ... \sum_{h_{2k-1} = 0}^{\lfloor \log\log X \rfloor} (\prod_{j=1}^{2k-1} \frac{\log X}{e^{\widetilde{h_j}}} \alpha(\widetilde{h_j}))^2 \int_{I_{\sigma,\textbf{h},m}} |F(\frac{1}{2} + it_{\sigma(2k)})|^2 \prod_{j=1}^{2k-1} \frac{|F(1/2+it_{\sigma(j)})|^2}{|F_{\widetilde{h_j}}(1/2+it_{\sigma(j)})|^2} dt_1 ... dt_{2k} , $$
where $\widetilde{h_j} := \min\{h_j, \lfloor \log\log X - 2\log\log\log X \rfloor\}$, say. (The introduction of $\widetilde{h_j}$ in place of $h_j$ is a technical manoeuvre simply to ensure that we always have $X^{e^{-\widetilde{h_j}}} \geq e^{(\log\log X)^2}$, which in particular is $\geq 1000k^6$, so we can efficiently control some error terms in the forthcoming calculation.)

Thus far all of our calculations have been deterministic, but now we shall exploit the randomness of the Euler products. For definiteness we work first in the Steinhaus case. We have
$$ |F(\frac{1}{2} + it_{\sigma(2k)})|^2 \prod_{j=1}^{2k-1} \frac{|F(\frac{1}{2}+it_{\sigma(j)})|^2}{|F_{\widetilde{h_j}}(\frac{1}{2}+it_{\sigma(j)})|^2} = \prod_{p \leq X} |1 - \frac{f(p)}{p^{1/2 + it_{\sigma(2k)}}}|^{-2} \prod_{j=1}^{2k-1} \prod_{X^{e^{-\widetilde{h_j}}} < p \leq X} |1 - \frac{f(p)}{p^{1/2 + it_{\sigma(j)}}}|^{-2} , $$
and here the primes smaller than $1000k^6$ only feature in the product corresponding to $t_{\sigma(2k)}$. So using the independence of the $f(p)$ for distinct primes $p$, together with the fact that (by the second part of Euler Product Result 1) $\E \prod_{p \leq 1000k^6} |1 - \frac{f(p)}{p^{1/2 + it_{\sigma(2k)}}}|^{-2} \ll \log(1000k^6)$, we deduce that the expectation of the products is
$$ \ll \log(1000k^6) \E \prod_{1000k^6 < p \leq X} |1 - \frac{f(p)}{p^{1/2 + it_{\sigma(2k)}}}|^{-2} \prod_{j=1}^{2k-1} \prod_{X^{e^{-\widetilde{h_j}}} < p \leq X} |1 - \frac{f(p)}{p^{1/2 + it_{\sigma(j)}}}|^{-2} . $$
Now we shall apply Euler Product Result 1 (and the independence of different $f(p)$) again to bound the expectation, this time exploiting the full uniformity of the result and using the fact that all the primes involved are $> 1000k^6$ to control the ``big Oh'' terms there. We obtain that the overall expectation is
$$ \ll \log(1000k^6) \exp\Biggl\{\sum_{j=1}^{2k-1} \sum_{\substack{X^{e^{-\widetilde{h_j}}} \\ < p \leq X}} \frac{1}{p} + \sum_{\substack{1000k^6 \\ < p \leq X}} \frac{1}{p} + 2\sum_{1 \leq j < l \leq 2k} \sum_{\substack{\max\{X^{e^{-\widetilde{h_j}}}, X^{e^{-\widetilde{h_l}}}\} \\ < p \leq X}} \frac{\cos((t_{\sigma(l)} - t_{\sigma(j)})\log p )}{p} \Biggr\} , $$
with the convention that when $l=2k$ the term $X^{e^{-\widetilde{h_{2k}}}}$ is omitted from the maxima. Note that we can also omit all terms $j,l$ from the double sum for which $h_j = 0$ or $h_l = 0$, since for those the inner sum over $p$ is empty. Using both parts of Number Theory Result 2, noting that all the differences $t_{\sigma(l)} - t_{\sigma(j)}$ here are $\ll 1$, we can further bound this expectation by
$$ \ll e^{O(k)} \log X \prod_{j=1}^{2k-1} e^{\widetilde{h_j}} \cdot \exp\Biggl\{7 \sum_{\substack{1 \leq j < l \leq 2k, \\ h_j, h_l \neq 0}} \frac{1}{(t_{\sigma(l)} - t_{\sigma(j)}) \max\{e^{-\widetilde{h_j}} \log X, e^{-\widetilde{h_l}} \log X\}}    \Biggr\} , $$
say, with the convention that $e^{-\widetilde{h_{2k}}} \log X$ is omitted from the maxima when $l=2k$.

Continuing with our analysis, using the facts that $\max\{e^{-\widetilde{h_j}} \log X, e^{-\widetilde{h_l}} \log X\} \geq \max\{e^{-h_j} \log X, e^{-h_l} \log X\}$ and that $\textbf{t} \in I_{\sigma,\textbf{h},m}$ (so whenever $h_i \neq 0$ we have $t_{\sigma(i+1)} - t_{\sigma(i)} \geq \frac{e^{h_i}}{\log X}$), we get a bound
$$ \ll e^{O(k)} \log X \prod_{j=1}^{2k-1} e^{\widetilde{h_j}} \cdot \exp\{7 \sum_{\substack{1 \leq j < l \leq 2k, \\ h_j, h_l \neq 0}} \frac{\min\{e^{h_j}, e^{h_l}\}}{e^{h_j} + \textbf{1}_{h_{j+1} \neq 0} e^{h_{j+1}} + ... + \textbf{1}_{h_{l-1} \neq 0} e^{h_{l-1}}} \} . $$
In the double sum, for fixed $j$ we can break the sum over those $j+1 \leq l \leq 2k$ with $h_l \neq 0$ into blocks, each containing as few consecutive terms as possible such that the sum of $e^{h_l}$ over each block is at least $e^{h_j}$. Then the total contribution to the $l$-sum from the $i$-th block of terms is $\leq \frac{2e^{h_j}}{i e^{h_j}} \leq 2/i$, so (since there are certainly at most $2k$ blocks) the sum over $l$ is $\leq 2\log(2k) + O(1)$, and the full double sum is $\leq 4k\log(2k) + O(k)$.

Note that although we wrote out the foregoing expectation calculations in the Steinhaus case and for $N=0$, things proceed exactly similarly in the Steinhaus case for general $N$. In the Rademacher case, using Euler Product Result 2 in place of Euler Product Result 1 one obtains the same terms along with an additional sum $2\sum_{\substack{1 \leq j < l \leq 2k, \\ h_j, h_l \neq 0}} \sum_{\max\{X^{e^{-\widetilde{h_j}}}, X^{e^{-\widetilde{h_l}}}\} < p \leq X} \frac{\cos((t_{\sigma(l)} + t_{\sigma(j)})\log p )}{p}$ in the exponential. Applying Number Theory Result 2 to this, recalling that we always have $|t_{\sigma(l)} + t_{\sigma(j)}| \ll |N| + 1 \ll (\log\log X)^2$ and $X^{e^{-\widetilde{h_j}}} \geq e^{(\log\log X)^2}$, this additional sum may be bounded by $7\sum_{\substack{1 \leq j < l \leq 2k, \\ h_j, h_l \neq 0}} \frac{1}{|t_{\sigma(l)} + t_{\sigma(j)}| \max\{e^{-\widetilde{h_j}} \log X, e^{-\widetilde{h_l}} \log X\}}$. And since all of the $t$ variables have the same sign (note that to ensure this when $N=0$, we made a little use of the fact that $|t|, |v| \geq (\log\log X)^{-2}$ in the definition of $\Delta_0$), we always have $|t_{\sigma(l)} + t_{\sigma(j)}| \geq t_{\sigma(l)} - t_{\sigma(j)}$ and this can be estimated in the same way as the existing sum in the exponential.

\vspace{12pt}
Putting things together, we have shown that the expectation of \eqref{integralsubeq} is
\begin{eqnarray}
& \ll & \log X \left(\frac{C\log\log\log X}{\log X}\right)^{2k} \sum_{|m| \leq k} \sum_{\sigma \in \text{Sym}(2k)} \nonumber \\
&& \sum_{h_1 = 0}^{\lfloor \log\log X \rfloor} ... \sum_{h_{2k-1} = 0}^{\lfloor \log\log X \rfloor} (\prod_{j=1}^{2k-1} \frac{\log X}{e^{\widetilde{h_j}}} \alpha(\widetilde{h_j}) )^2 \int_{I_{\sigma,\textbf{h},m}} e^{56k\log(2k) + O(k)} \log X \prod_{j=1}^{2k-1} e^{\widetilde{h_j}} dt_1 ... dt_{2k} \nonumber \\
& = & \log X (C k^{28} \log\log\log X)^{2k} \sum_{h_1, ..., h_{2k-1} = 0}^{\lfloor \log\log X \rfloor} (\prod_{j=1}^{2k-1} \frac{\log X}{e^{\widetilde{h_j}}} \alpha(\widetilde{h_j})^{2} ) \sum_{|m| \leq k} \sum_{\sigma \in \text{Sym}(2k)} \int_{I_{\sigma,\textbf{h},m}} dt_1 ... dt_{2k} . \nonumber
\end{eqnarray}
When estimating the inner sums here, recalling the definition of $I_{\sigma,\textbf{h},m} \subseteq I_m$ we always have the easy bound $\int_{I_{\sigma,\textbf{h},m}} dt_1 ... dt_{2k} \leq \prod_{j=1}^{2k-1} (\frac{e^{h_j + 1}}{\log X})$. Alternatively, for any given $1 \leq l \leq 2k$ we can upper bound $\int_{I_{\sigma,\textbf{h},m}} dt_1 ... dt_{2k}$ by
$$ \text{meas}\{\textbf{t} \in [0,1]^{2k} : t_{\sigma(j+1)} - t_{\sigma(j)} \leq \frac{e^{h_j + 1}}{\log X} \; \forall \; j \neq l-1, l, \;\;\; |(t_1 + ... + t_k - t_{k+1} - ... - t_{2k}) - m| \leq \frac{2}{\log^{2/3}X} \} , $$
which (regarding the final inequality as a condition on $t_{\sigma(l)}$) is all $\ll \frac{1}{\log^{2/3}X} \prod_{\substack{j=1, \\ j \neq l-1, l}}^{2k-1} (\frac{e^{h_j + 1}}{\log X})$. Combining these bounds, noting that summing over $m$ and $\sigma$ produces another factor $k(2k)! \ll (2k)^{2k}$ and that $\frac{e^{h_j}}{e^{\widetilde{h_j}}} \leq e(\log\log X)^2$ for all $h_j$ as well, we deduce that the expectation of \eqref{integralsubeq} is
$$ \ll \log X (C k^{29} (\log\log X)^2 \log\log\log X)^{2k} \sum_{h_1, ..., h_{2k-1} = 0}^{\lfloor \log\log X \rfloor} (\prod_{j=1}^{2k-1} \alpha(\widetilde{h_j})^{2} ) \min_{1 \leq l \leq 2k}\{1, \frac{1}{\log^{2/3}X} \frac{\log^{2}X}{e^{h_{l-1}+h_l}} \} . $$
In the minimum here, we apply the convention that when $l=1$ or $l=2k$, the term $\frac{1}{\log^{2/3}X} \frac{\log^{2}X}{e^{h_{l-1}+h_l}}$ is replaced by $\frac{1}{\log^{2/3}X} \frac{\log X}{e^{h_1}}$ and $\frac{1}{\log^{2/3}X} \frac{\log X}{e^{h_{2k-1}}}$ respectively.

Finally, note that if $\textbf{h}$ is a tuple in the multiple sum for which there exists some $1 \leq l \leq 2k$ with $\frac{\log X}{e^{h_{l-1}}} \leq \log^{1/6}X$ and $\frac{\log X}{e^{h_{l}}} \leq \log^{1/6}X$, then $\min_{1 \leq l \leq 2k}\{1, \frac{1}{\log^{2/3}X} \frac{\log^{2}X}{e^{h_{l-1}+h_l}} \} \leq \frac{1}{\log^{1/3}X}$. In this case we can just apply the trivial bound $\prod_{j=1}^{2k-1} \alpha(\widetilde{h_j})^{2} \leq (\log\log X)^{24k}$. Alternatively, if $\textbf{h}$ is a tuple for which no such $l$ exists then for at least half of all $1 \leq l \leq 2k-1$ we certainly have $\widetilde{h_l} \leq h_l \leq 0.99\log\log X$, and so $\prod_{j=1}^{2k-1} \alpha(\widetilde{h_j})^{2} \leq (\frac{1}{(\log\log X)^{2000}})^{k} (\log\log X)^{12(k-1)} \leq (\frac{1}{(\log\log X)^{994}})^{2k}$. Overall, noting also that the total number of tuples $\textbf{h}$ is $(\lfloor \log\log X \rfloor + 1)^{2k-1} \leq (2\log\log X)^{2k}$, we deduce that the expectation of \eqref{integralsubeq} is
\begin{equation}\label{finalexpeq}
\ll \log X \left(\frac{1}{\log^{1/3}X} (C k^{29} (\log\log X)^{15} \log\log\log X)^{2k} + (\frac{C k^{29} \log\log\log X}{(\log\log X)^{991}})^{2k} \right) .
\end{equation}
One obtains a comparable bound for all $-(\log\log X)^2 - 1 \leq N \leq (\log\log X)^2$, and so a comparable bound when summing over $N$ with the prefactor $\frac{1}{N^2 + 1}$. So by Markov's inequality, with probability $\geq 1 - (\log\log X)^{-2k}$ the part of our original integrals where $|| t_1 + ... + t_k - t_{k+1} - ... - t_{2k} || < 1/\log^{2/3}X$ will be at most $(\log\log X)^{2k}$ multiplied by \eqref{finalexpeq}, which (after simplifying the fractions a little) finishes the proof.
\end{proof}

%%\vspace{12pt}
%%\noindent {\em Acknowledgements.}


\begin{thebibliography}{99}

\bibitem{basquin} J. Basquin. Sommes friables de fonctions multiplicatives al\'eatoires. {\em Acta Arith.}, \textbf{152}, no. 3, pp 243-266. 2012

\bibitem{bobgoldgrankou} J. Bober, L. Goldmakher, A. Granville, D. Koukoulopoulos. The frequency and the structure of large character sums. {\em J. Eur. Math. Soc. (JEMS)}, \textbf{20}, no. 7, pp 1759-1818. 2018

\bibitem{chandee} V. Chandee. On the correlation of shifted values of the Riemann zeta function. {\em Q. J. Math.}, \textbf{62},  no. 3, pp 545-572. 2011

\bibitem{erdos} P. Erd\H{o}s. Some applications of probability methods to number theory. In {\em Proc. of the 4th Pannonian Symp. on Math. Stat., (Bad Tatzmannsdorf, Austria 1983)}. {\em Mathematical statistics and applications}, vol. B, pp 1-18. Reidel, Dordrecht. 1985

\bibitem{gransoundlcs} A. Granville, K. Soundararajan. Large character sums. {\em J. Amer. Math. Soc.}, \textbf{14}, no. 2, pp 365-397. 2001

\bibitem{gut} A. Gut. {\em Probability: A Graduate Course.} Second edition, published by Springer Texts in Statistics. 2013

\bibitem{halasz} G. Hal\'{a}sz. On random multiplicative functions. In {\em Hubert Delange Colloquium, (Orsay, 1982)}. {\em Publications Math\'{e}matiques d'Orsay}, \textbf{83}, pp 74-96. Univ. Paris XI, Orsay. 1983

\bibitem{harpergp} A. J. Harper. Bounds on the suprema of Gaussian processes, and omega results for the sum of a random multiplicative function. {\em Ann. Appl. Probab.}, \textbf{23}, no. 2, pp 584-616. 2013

\bibitem{harperlcz} A. J. Harper. A note on the maximum of the Riemann zeta function, and log-correlated random variables. Preprint available online at \url{https://arxiv.org/abs/1304.0677}

\bibitem{harperrmfhigh} A. J. Harper. Moments of random multiplicative functions, II: High moments. {\em Algebra Number Theory}, \textbf{13}, no. 10, pp 2277-2321. 2019

\bibitem{harperpartition} A. J. Harper. On the partition function of the Riemann zeta function, and the Fyodorov--Hiary--Keating conjecture. Preprint available online at \url{https://arxiv.org/abs/1906.05783}

\bibitem{harperrmflow} A. J. Harper. Moments of random multiplicative functions, I: Low moments, better than squareroot cancellation, and critical multiplicative chaos. {\em Forum of Mathematics, Pi}, \textbf{8}, e1, 95pp. 2020

\bibitem{harperlamzouri} A. J. Harper, Y. Lamzouri. Orderings of weakly correlated random variables, and prime number races with many contestants. {\em Probab. Theory Related Fields}, \textbf{170}, no. 3-4, pp 961-1010. 2018

\bibitem{tenenbaum} Y-K. Lau, G. Tenenbaum, J. Wu. On mean values of random multiplicative functions. {\em Proc. Amer. Math. Soc.}, \textbf{141}, pp 409-420. 2013. Also see \url{www.iecl.univ-lorraine.fr/~Gerald.Tenenbaum/PUBLIC/Prepublications_et_publications/RMF.pdf} for some corrections to the published version.

\bibitem{mont} H. L. Montgomery. {\em Ten Lectures on the Interface Between Analytic Number Theory and Harmonic Analysis.} Published for the Conference Board of the Mathematical Sciences by the American Mathematical Society. 1994

\bibitem{mv} H. L. Montgomery, R. C. Vaughan. {\em Multiplicative Number Theory I: Classical Theory.} First edition, published by Cambridge University Press. 2007

\bibitem{ng} N. Ng. The distribution of the summatory function of the M\"{o}bius function. {\em Proc. London. Math. Soc.}, \textbf{89}, no. 3, pp. 361-389, 2004.

\bibitem{rr} G. Reinert, A. R\"{o}llin. Multivariate Normal Approximation with Stein's Method of Exchangeable Pairs under a General Linearity Condition. {\em Annals of Probability}, \textbf{37}, no. 6, pp 2150-2173. 2009

\bibitem{soundyoung} K. Soundararajan, M. P. Young. The second moment of quadratic twists of modular $L$-functions. {\em J. Eur. Math. Soc. (JEMS)}, \textbf{12}, no. 5, pp 1097-1116. 2010

\bibitem{tenenbaummultmean} G. Tenenbaum. Moyennes effectives de fonctions multiplicatives complexes. {\em Ramanujan J.}, \textbf{44}, no. 3, pp 641-701. 2017

\bibitem{wintner} A. Wintner. Random factorizations and Riemann's hypothesis. {\em Duke Math. J.}, \textbf{11}, pp 267-275. 1944



%%\bibitem{bbsszpseudo} A. Bondarenko, O. F. Brevig, E. Saksman, K. Seip, J. Zhao. Pseudomoments of the Riemann zeta function. Preprint available online at {\verb https://arxiv.org/abs/1701.06842v2 }.











%%\bibitem{saksmanseip} E. Saksman, K. Seip. Some open questions in analysis for Dirichlet series. In {\em Recent progress on operator theory and approximation in spaces of analytic functions}, pp 179-191, Contemp. Math., 679, Amer. Math. Soc., Providence, RI. 2016




\end{thebibliography}
\end{document}